\documentclass[11pt]{amsart}
\usepackage{amsmath}
\usepackage{amssymb}
\usepackage{latexsym}
\usepackage{amsfonts}

\usepackage{xcolor}

\usepackage[plainpages=false, hypertexnames=false, pdfpagelabels=true, hyperindex=true, linktocpage, pagebackref=true, pdfa=true]{hyperref}

\usepackage{amsrefs}

\begin{document}
\newcommand{\ci}[1]{_{ {}_{\scriptstyle #1}}}
\newcommand{\norm}[1]{\ensuremath{\|#1\|}}
\newcommand{\abs}[1]{\ensuremath{\vert#1\vert}}
\newcommand{\p}{\ensuremath{\partial}}
\newcommand{\pbar}{\ensuremath{\bar{\partial}}}
\newcommand{\db}{\overline\partial}
\newcommand{\D}{\mathbb{D}}
\newcommand{\T}{\mathbb{T}}
\newcommand{\C}{\mathbb{C}}
\newcommand{\N}{\mathbb{N}}
\newcommand{\td}{\widetilde\Delta}
\newcommand{\La}{\langle }
\newcommand{\Ra}{\rangle }
\newcommand{\tr}{\operatorname{tr}}
\newcommand{\ran}{\operatorname{Ran}}
\newcommand{\vf}{\varphi}
\newcommand{\e}{\varepsilon}
\newcommand{\f}[2]{\ensuremath{\frac{#1}{#2}}}
\newcommand{\be}{\mathbf{e}}
\newcommand{\clos}{\operatorname{clos}}
\newcommand{\rank}{\operatorname{rank}}
\newcommand{\bz}{\mathbf{z}}
\newcommand{\tto}{\!\!\to\!}
\newcommand{\wt}{\widetilde}
\newcommand{\shto}{\raisebox{.3ex}{$\scriptscriptstyle\rightarrow$}\!}
\newcommand{\bL}{\mathbf{L}}
\newcommand{\entrylabel}[1]{\mbox{#1}\hfill}




\numberwithin{equation}{section}
\newtheorem{thm}{Theorem}[section]
\newtheorem{lm}[thm]{Lemma}
\newtheorem{cor}[thm]{Corollary}
\newtheorem{prop}[thm]{Proposition}
\newtheorem{rem}[thm]{Remark}
\newtheorem*{rem*}{Remark}

\title{The Corona Problem for Slice Hyperholomorphic Functions}

\author[F. Colombo]{Fabrizio Colombo}
\address{(FC)
Politecnico di Milano\\Dipartimento di Matematica\\Via E. Bonardi 9\\20133
Milano, Italy}
\email{fabrizio.colombo@polimi.it}

\author[E. Pozzi]{Elodie Pozzi}
\address{(EP) Department of Mathematics and Statistics \\
Saint Louis University \\
220 N. Grand Blvd, 63103 St Louis MO, USA}
\email{elodie.pozzi@slu.edu}

\author[I. Sabadini]{Irene Sabadini}
\address{(IS)
Politecnico di Milano\\Dipartimento di Matematica\\Via E. Bonardi 9\\20133
Milano, Italy
}
\email{irene.sabadini@polimi.it}

\author[B. D Wick]{Brett D. Wick}
\address{(BDW)  Department of Mathematics, Washington University -
 St. Louis, One Brookings Drive\\ St. Louis, MO USA 63130-4899}
\email{wick@math.wustl.edu}
\thanks{BDW's research supported in part by National Science Foundation DMS awards \#2349868 and \#2054863 and Australian Research Council -- DP 220100285.}

\thanks{EP and BDW would like to thank Politecnico di Milano for hospitality during visits in Spring 2024 when this project was started.}

\begin{abstract}
This paper addresses the Corona problem for slice hyperholomorphic functions for a single quaternionic variable. While the Corona problem is well-understood in the context of one complex variable, it remains highly challenging in the case of several complex variables. 
The extension of the theory of one complex variable to several complex variables is not the only possible extension to multi-dimensional complex analysis.  Instead of functions holomorphic in each variable separately, in this paper we will consider functions in some hypercomplex algebras, in particular in the algebra of quaternions.
\\\indent Previously, the Corona problem had not been studied within the hypercomplex framework because of challenges posed by pointwise multiplication, which is not closed for hypercomplex-valued analytic functions. Alternative notions of multiplication that are closed often compromise other desirable properties, further complicating the analysis.
\\\indent In this work, we resolve the Corona problem within the quaternionic slice hyperholomorphic setting. Our approach involves reformulating the quaternionic Bezout equation with respect to the appropriate multiplication into a new system of Bezout equations on the unit disc. We solve this system by adapting Wolff's proof of the Corona theorem for bounded analytic functions. As the number of generators increases, the associated algebra grows increasingly intricate.
\end{abstract}

\maketitle



	\tableofcontents

\section{Introduction and Statement of Main Results}
\label{s:IntroMain}

Carleson's Corona Theorem for the algebra of bounded analytic functions on the unit disk, $H^\infty(\mathbb{D})$, has proven to be a cornerstone in complex analysis because of its connections to function theory, operator theory, and harmonic analysis.  For a collection $f_1,\ldots, f_n\in H^\infty(\mathbb{D})$, the Corona problem, is to understand the equivalence between the Bezout equation and the so called Carleson condition:
\begin{align*}
    \sum_{j=1}^n f_j(z)g_j(z) &=1\quad\forall z\in \mathbb{D} \quad\textnormal{ and }\quad 0<\delta^2\leq\sum_{j=1}^{n}\vert f_j(z)\vert^2\leq 1\quad z\in\mathbb{D}.
\end{align*}
If the Bezout equation holds then it must be the case that the $f_j$ do not simultaneously vanish, which is quantified by the second condition above. Surprisingly Carleson demonstrated that the functions not simultaneously vanishing was sufficient to solve the Bezout equation for solutions $g_j\in H^\infty(\mathbb{D})$, see \cite{MR0141789}.  The formulation of this question reveals its connections to analysis, algebra, and topology, offering a diverse and rich array of tools and techniques for tackling problems in this area; see the conference proceedings \cite{MR3329539} devoted to this question and some of its aspects.

While much progress has been made in studying the Corona problem in one complex variable, the related question in several complex variables has proven to be more intractable.  
In the context of several complex variables, there exist domains where the $H^\infty$ Corona problem is known to fail, \cite{MR0916722},  and there are very simple domains such as the unit ball or the polydisc for which the Corona problem remains open.
 Variants and generalizations of this question have been extensively well studied.  We point to the bibliography for a sampling of generalizations on the Corona problem.

\medskip
The extension of the theory of one complex variable to several complex variables is not the only possible multi-dimensional framework.
In fact, there exists an extension procedure called the Fueter-Sce-Qian mapping theorem, see the book \cite{MR4240465} and the references therein,
 that leads to two different notions of hyperholomorphic functions: slice hyperholomorphic functions
 see \cites{MR3585395,MR2752913,MR3013643},
  and monogenic functions, see
 \cites{MR2089988,MR3496962,MR1169463}.
Monogenic functions have been investigated since the beginning of the last century,
while the theory of slice hyperholomorphic functions started in 2006. These function theories are both very well studied and find numerous applications in operator theory and harmonic analysis. In fact,
 slice hyperholomorphic functions generate a spectral theory based on the $S$-spectrum, see \cites{MR3967697,MR3887616,MR2752913},
while monogenic functions induce the spectral theory based on the monogenic spectrum, see \cites{MR3930595,MR2069293}.

There is also an interesting relationship between the two hyperholomorphic spectral theories via the $F$-functional calculus and the recent theory of the fine structures on the $S$-spectrum has given
a further impulse in the expansion of both the function theories, which now include also polyanalytic and polyharmonic functions, and the spectral theories, see \cites{MR4502763,MR4576314,MR4674273,MR4741057}.
Moreover, the two main hyperholomorphic spectral theories have different and complementary applications.

Several results available for holomorphic functions in complex analysis are also valid in the hypercomplex setting, possibly with variations in the statements and/or in the techniques of their proofs. Questions pertaining the study of ideals of hyperholomorphic functions and the Corona problems are intrinsically problematic and especially challenging since the pointwise multiplication of functions does not preserve hyperholomorphicity. Suitable products can be introduced in order to have a closed operation in the class of hyperholomorphic functions, but not all the properties that one expects for a product are still valid. 

\medskip
In this paper we prove a Corona theorem for slice hyperholomorphic functions on the unit ball of the quaternions.  We can solve the Corona problem in this context because some of the techniques and tools from complex analysis are more readily amenable to adaption in the quaternionic setting.  An earlier, related, version of the result is given in \cite{MR3621101} where the authors studied ideals in the quaternionic setting.  They assumed that the functions appearing in the Bezout equation do not have common zeros, but were not able to obtain norm estimates on the resulting solutions because of the use of sheaf theory. In this paper we consider instead a quantitative version of the problem, assuming the Carleson condition on the Corona data.

This Corona theorem in the slice hyperholomorphic case will follow as a corollary to the following new Bezout-type Corona theorem on the unit disk.  To state the result, we let $\hat{P}(z)=\overline{P(\overline{z})}$, noting that $\hat{\cdot}$ acts as an involution on the algebra $H^\infty(\mathbb{D})$.  In particular, if $P$ is holomorphic, then so is the function $\hat{P}(z)$, and if $P$ is bounded, then so is $\hat{P}$ with the same bound.  

Our new Corona theorem is as follows:
\begin{thm}
    \label{t:NewCorona}
    Let $0<\delta<1$ and $n\in\mathbb{N}$ and $\mathbb{D}$ be the unit disk in the complex plane.  Suppose that $F_1,\ldots, F_n$ and $G_1,\ldots, G_n\in H^\infty(\mathbb{D})$ and satisfy
\begin{align}
\label{e:NewCC}
0<\delta^2 & \leq \sum_{r=1}^{n}\sum_{j=1}^n
\left\vert F_j(z)\hat{F}_r(z)+G_r(z)\hat{G}_j(z)\right\vert^2
+\sum_{r=1}^{n}\sum_{j=r+1}^n\left\vert \hat{G}_r(z)F_j(z)-\hat{G}_j(z)F_r(z)\right\vert^2\\
& +\sum_{r=1}^{n}\sum_{j=r+1}^n\left\vert \hat{F}_r(z)G_j(z)-\hat{F}_j(z)G_r(z)\right\vert^2\leq 1\quad\forall z\in\mathbb{D}\notag
\end{align}
    and
    \begin{equation}
    \label{e:OldCC}
        0<\delta^2\leq\sum_{j=1}^n\vert F_j(z)\vert^2+\sum_{j=1}^n\vert G_j(z)\vert^2\leq 1\quad\forall z\in\mathbb{D}.
    \end{equation}
    Then there exist functions $H_1,\ldots, H_n$ and $K_1,\ldots, K_n\in H^\infty(\mathbb{D})$ so that
    \begin{align*}
    \sum_{j=1}^n F_jH_j-\sum_{j=1}^nG_j\hat{K}_j & =1\\
    \sum_{j=1}^n F_jK_j+\sum_{j=1}^nG_j\hat{H}_j & =0,
\end{align*}
and
$$\displaystyle\max_{1\leq j\leq n}\left\{\Vert H_j\Vert_{H^\infty(\mathbb{D})},\Vert K_j\Vert_{H^\infty(\mathbb{D})}\right\}\lesssim C(\delta,n):= \frac{1}{\delta^2}+\frac{n}{\delta^4}+\left(\frac{1}{\delta^2}
+\frac{n}{\delta^4}\right)\left(\frac{n}{\delta^4}+\frac{n^2}{\delta^8}\right).
$$
\end{thm}

Conditions like \eqref{e:NewCC} and \eqref{e:OldCC} are the standard ``Corona conditions'' that indicate the functions do not simultaneously vanish.  Condition \eqref{e:OldCC} is the well-known standard condition one would encounter in the solution to these types of problems.  While condition \eqref{e:NewCC} is what rules out potential zeros that could lead to the systems of equations having no solution; these are motivated additionally by the zero condition for slice regular functions with a quaternionic variable.  Both condition are in fact necessary, \eqref{e:OldCC} follows from a standard argument, while \eqref{e:NewCC} will be seen as necessary in the course of the argument.  Further, since $n$ is finite, then applications of Cauchy-Schwarz and equivalence of $\ell^p$ norms on finite dimensional spaces show that \eqref{e:NewCC} implies \eqref{e:OldCC} with $\delta'=\delta'(\delta,n)$.  One could additionally assume that conditions \eqref{e:NewCC} and \eqref{e:OldCC} hold with values $\delta_1$ and $\delta_2$ and then deduce an estimate on $C(\delta_1,\delta_2,n)$ from these conditions; we do not presently explore this.\\

As an application of this result, in Theorem \ref{t:QCorona} we demonstrate the Corona theorem for slice hyperholomorphic functions in the quaternionic setting.  This resolves an open  question in the area of quaternionic analysis, and a question specifically raised in \cite{Saracco}; related results under a different non-vanishing condition can be seen in \cite{Shelah_MSThesis}.

As already discussed, the product of two slice hyperholomorphic functions is not, in general, slice hyperholomorphic. In order to preserve slice hyperholomorphicity, a suitable product denoted by $\ast$ must be defined. To avoid technicalities, we take the following as definition:
$$
(f\ast g)(q):=f(q)g(f^{-1}(q)qf(q))
$$
when $f(q)\neq 0$ and zero when $f(q)=0$.  Because the product $\ast$ is not commutative one must decide whether to multiply on the left or right in the Bezout equation.  Already in \cite{Shelah_MSThesis} it was observed that seeking right inverses in the Bezout equation is the correct question.

Let $H^\infty(\mathbb{B})$ denote the algebra of slice hyperholomorphic functions on the unit ball $\mathbb{B}$ in the quaternions.  It is normed by the usual supremum norm:
$$
\left\Vert f\right\Vert_{H^{\infty}(\mathbb{B})}=\sup_{q\in\mathbb{B}}\vert f(q)\vert.
$$
See Section \ref{s:SlicePrelims} for more details and background on this area of complex analysis.  Our quantitative version of the Corona result in the slice hyperholomorphic setting is the following:
\begin{thm}
\label{t:QCorona}
    Let $0<\delta<1$ and $n\in\mathbb{N}$ and let $\mathbb{B}$ be the unit ball in the quaternions.  Suppose that $f_1,\ldots, f_n\in H^\infty(\mathbb{B})$ are slice hyperholomorphic functions such that
    $$
        0<\delta^2\leq \sum_{j=1}^n\left\vert f_j(q)\right\vert^2\leq 1\quad\forall q\in\mathbb{B}.
    $$
Then there are $g_1,\ldots, g_n$ such that
    $$
    \sum_{j=1}^n (f_j \ast g_j) (q)=1\quad\forall q\in\mathbb{B}
    $$
    and $$\displaystyle\max_{1\leq j\leq n} \Vert g_j\Vert_{H^\infty(\mathbb{B})}\lesssim C(\delta,n):= \frac{1}{\delta^2}+\frac{n}{\delta^4}+\left(\frac{1}{\delta^2}
    +\frac{n}{\delta^4}\right)\left(\frac{n}{\delta^4}+\frac{n^2}{\delta^8}\right).
    $$
\end{thm}
This is an example of solvability of the Corona problem in a $4$-dimensional case which has the potential for an extension to the $n$-dimensional case for Clifford algebra valued functions.  See some of the questions raised in the last section of the paper for more details.

\medskip
In Section \ref{s:NewCorona}, we introduce the new Corona problem which is inspired by the splitting of quaternions into real and imaginary parts, and whose solution leads to Theorem \ref{t:NewCorona}. In the first subsection, we provide an explanation and proof in the case of two generators as the algebra is easier in this setting. The next subsection addresses the case of more generators, and contains additional technicalities. In the last part we show how to perturb the smooth and bounded solutions to the problem in order to obtain holomorphic solutions and we adapt and modify Wolff's proof of the Corona theorem in this new setting.

Section \ref{s:QCorona} contains the preliminaries on slice hyperholomorphic functions and the statement of the $H^\infty(\mathbb{B})$-Corona problem. Then we prove the result first in the case of one generator and then in the general case, namely Theorem \ref{t:QCorona}.

Finally, in Section \ref{s:Conclude} we state some concluding remarks and potential lines of further investigation.

\section{A New Corona Problem}
\label{s:NewCorona}

In this section we prove the following new Corona theorem, restated from the introduction.

\begin{thm}
    \label{t:NewCoronaS2}
    Let $0<\delta<1$ and $n\in\mathbb{N}$.  Suppose that $F_1,\ldots, F_n$ and $G_1,\ldots, G_n\in H^\infty(\mathbb{D})$ and satisfy
\begin{align}
\label{e:NewCC2}
0<\delta^2 & \leq \sum_{r=1}^{n}\sum_{j=1}^n
\left\vert F_j(z)\hat{F}_r(z)+G_r(z)\hat{G}_j(z)\right\vert^2
+\sum_{r=1}^{n}\sum_{j=r+1}^n\left\vert \hat{G}_r(z)F_j(z)-\hat{G}_j(z)F_r(z)\right\vert^2\\
& +\sum_{r=1}^{n}\sum_{j=r+1}^n\left\vert \hat{F}_r(z)G_j(z)-\hat{F}_j(z)G_r(z)\right\vert^2\leq 1\quad\forall z\in\mathbb{D}\notag
\end{align}
and
    \begin{equation}
    \label{e:CoronaCondition}
        0<\delta^2\leq\sum_{j=1}^n\vert F_j(z)\vert^2+\sum_{j=1}^n\vert G_j(z)\vert^2\leq 1\quad\forall z\in\mathbb{D}.
    \end{equation}
    Then there exist functions $H_1,\ldots, H_n$ and $K_1,\ldots, K_n\in H^\infty(\mathbb{D})$ so that
    \begin{align}
    \sum_{j=1}^n F_jH_j-\sum_{j=1}^nG_j\hat{K}_j & =1\label{e:nNewBezout11}\\
    \sum_{j=1}^n F_jK_j+\sum_{j=1}^nG_j\hat{H}_j & =0\label{e:nNewBezout22}
\end{align}
and
\begin{equation}\label{e:nNewBezout33}\displaystyle\max_{1\leq j\leq n}\left\{\Vert H_j\Vert_{H^\infty(\mathbb{D})},\Vert K_j\Vert_{H^\infty(\mathbb{D})}\right\}\lesssim C(\delta,n):= \frac{1}{\delta^2}+\frac{n}{\delta^4}+\left(\frac{1}{\delta^2}
+\frac{n}{\delta^4}\right)\left(\frac{n}{\delta^4}+\frac{n^2}{\delta^8}\right).
\end{equation}
\end{thm}
First of all, we observe that \eqref{e:CoronaCondition} follows from \eqref{e:nNewBezout11} via the standard argument.  We will argue in the course of the proof below that \eqref{e:nNewBezout11} and \eqref{e:nNewBezout22} imply that condition \eqref{e:NewCC2} is necessary as well.

The general idea for how to prove this result is to follow a similar strategy to prove the classical result of Carleson, but using idea's from Wolff's proof of the result.  These ideas can be found in the monographs \cites{MR1419088, MR0827223, MR1669574, MR2261424}.

More specifically, the idea will be to first construct smooth bounded solutions $h_j$ and $k_j$ to \eqref{e:nNewBezout11} and \eqref{e:nNewBezout22}.  For the equation \eqref{e:nNewBezout11} this will be immediate from \eqref{e:CoronaCondition}, while for the other it will be a result of some linear algebra arguments that it is possible to solve \eqref{e:nNewBezout11} and \eqref{e:nNewBezout22} simultaneously.  In particular, this is accomplished by looking at certain syzygies associated to the system of equations.  These smooth bounded solutions will then be perturbed to make the resulting solutions holomorphic.  This will be accomplished by setting up a system of $\overline{\partial}$-equations related to the smooth and bounded solutions.  In the classical case, these arise from the Koszul complex; while in this case the system of equations comes from the quaternionic setting, which is non-commutative. Thus the solutions have to be determined by hand and shown to solve the equations of interest.  The bounded and holomorphic solutions will then follow from Wolff's ideas for bounded solutions to $\overline{\partial}$-equations.  While the general approach follows the classical proof, difficulties emerge due to the additional algebraic complexities associated with the two equations \eqref{e:nNewBezout11} and \eqref{e:nNewBezout22}.  Miraculously, these two equations yield a new Bezout equation, \eqref{e:nNewBezout33}, that plays a key algebraic role in the solution, as well as additional identities that are pivotal to showing the holomorphicity of the solutions. Condition \eqref{e:nNewBezout33} will be a necessary condition, and it in turn will show that \eqref{e:NewCC2} must hold as well.

In what follows, we will let $h_j$ and $k_j$ denote bounded and smooth, {\it i.e.} $C^{\infty}$, solutions to equations \eqref{e:nNewBezout11} and \eqref{e:nNewBezout22}, and will let $H_j$ and $K_j$ denote bounded and holomorphic solutions to the equations \eqref{e:nNewBezout11} and \eqref{e:nNewBezout22}.

\subsection{Explanation and Proof in the Case of \texorpdfstring{$n=1$}{n is 1}}
The case of $n=1$ is addressed separately since it helps to illustrate some of the challenges in this question.  The arguments provided will be a little brief to show some of the main ideas and challenges, while in the general case full details will be carried out in Section \ref{s:n=2} and Section \ref{s:3+GeneratorsThm} where there is more algebraic challenges to address.  The point is to highlight the ideas so the reader can more easily follow the arguments later.

Since $n=1$, the equations to solve are:
\begin{align}
    F_1H_1-G_1\hat{K}_1 &=1 \label{e:eq1n=1}\\
    F_1K_1+G_1\hat{H}_1 & = 0\label{e:eq2n=1}.
\end{align}
The system of equations \eqref{e:eq1n=1} and \eqref{e:eq2n=1} imply another equation; namely any solutions of these systems must also satisfy the equation:
\begin{align}
    (\hat{F}_1F_1+\hat{G}_1G_1)(\hat{H}_1H_1+\hat{K}_1K_1)=1\label{e:eq3n=1}.
\end{align}
This follow because by simple algebra, since we can multiply equations \eqref{e:eq1n=1} and \eqref{e:eq2n=1} by $F$ or $G$ and using the involution $\hat{\cdot}$ we can arrive at the identities:
\begin{align}
    F_1 & = (\hat{F}_1F_1+\hat{G}G_1)\hat{H}_1\label{e:n=1Identity}\\
    G_1 & = -(\hat{F}_1F_1+\hat{G}_1G_1)\hat{K}_1.\label{e:n=1Identity2}
\end{align}
These identities can be substituted back into \eqref{e:eq1n=1} to arrive at \eqref{e:eq3n=1}.  This idea will appear when we handle the case of general $n \geq 2$, but the resulting identities become more algebraically involved.

Now observe that \eqref{e:eq3n=1} implies the condition:
\begin{equation}
\label{e:NewCCn=1}
0<\delta^2\leq \vert \hat{F}_1(z)F_1(z)+\hat{G}_1(z)G_1(z)\vert^2\leq 1\quad\forall z\in\mathbb{D}
\end{equation}
which is the analog of condition \eqref{e:NewCC2} in the case $n=1$ since the other terms in that equation vanish.  Define the function $\Delta(z):=\vert \hat{F}_1(z)F_1(z)+\hat{G}_1(z)G_1(z)\vert^2$,  and then it holds that $0<\delta^2\leq \Delta(z)\leq 1$.

The idea is now to build smooth solutions to equations \eqref{e:eq1n=1} and \eqref{e:eq2n=1}.  Define $D(z):=\left\vert F_1(z)\right\vert^2+\left\vert G_1(z)\right\vert^2$, by \eqref{e:CoronaCondition} we have that $0<\delta^2\leq D(z)\leq 1$.  Define the functions:
\begin{align*}
\varphi_1(z)&:=\frac{\overline{F_1(z)}}{D(z)} & \psi_1(z) &:=-\frac{G_1(\overline{z})}{D(\overline{z})}.
\end{align*}
Then the functions $\varphi$ and $\psi$ belong to $L^\infty(\mathbb{D})$ with norm controlled by a function of $\delta$.   It is also immediate via algebra that:
$$
F_1\varphi_1-G_1\hat{\psi}_1=1.
$$
Then perturb the solutions to create smooth and bounded solutions to the problem:
\begin{align}
    h_1(z) & := \varphi_1(z)+\lambda_1(z)G_1(z)\\
    k_1(z) & :=\psi_1(z)+\hat{\lambda}_1(z)\hat{F}_1(z).
\end{align}
It then is possible to show that \eqref{e:eq1n=1} holds:
$$
F_1h_1-G_1\hat{k}_1=1
$$
for any choice of smooth bounded function $\lambda_1$ with an additional condition to hold.  Now observe that it is possible to choose the function $\lambda_1$ to force \eqref{e:eq2n=1} to hold too.  Again via simple algebra \eqref{e:eq2n=1} is true, namely
$$
F_1k_1+G_1\hat{h}_1=0
$$
if and only if it is possible to solve
$$
F_1\psi_1+G_1\hat{\varphi}_1+\hat{\lambda}_1\left(\hat{F}_1F_1+\hat{G}_1G_1\right)=0.
$$
By \eqref{e:NewCCn=1} simply choose 
$$
\hat{\lambda}_1:=-\frac{F_1\psi_1+G_1\hat{\varphi}_1}{\left(\hat{F}_1F_1+\hat{G}_1G_1\right)}
$$
and this function will be smooth and bounded since $F_1$, $G_1$, $\varphi_1$ and $\psi_1$ are and \eqref{e:NewCCn=1} holds.

Finally we show these solutions are in fact holomorphic.  In the case of general $n$ this requires additional work and establishing that a system of $\overline{\partial}$-equations holds have solutions.

To accomplish this set
\begin{align*}
    H_1(z) & := h_1(z)\\
    K_1(z) & := k_1(z).
\end{align*}
Then by the construction above we have that $H_1$ and $K_1$ satisfy \eqref{e:eq1n=1} and \eqref{e:eq2n=1} and additionally are bounded by some constant depending only on $\delta$.  Surprisingly these solutions are also holomorphic.  To see this, we simply need to take $\overline{\partial}$ of each of these equations (after first taking $\hat{\cdot}$).  As they are similar we consider the first one and taking $\overline{\partial}$ gives, since $F_1, G_1, \hat{F}_1$ and $\hat{G}_1$ are holomorphic that:
$$
0=\left(\hat{F}_1F_1+\hat{G}_1G_1\right)\overline{\partial}H_1
$$
and since by \eqref{e:NewCCn=1} we have that $\left(\hat{F}_1F_1+\hat{G}_1G_1\right)$ is invertible we arrive at $\overline{\partial}H_1=0$ and so $H_1$ is holomorphic as well.  For the function $K_1$ the argument is identical.  These are then the bounded analytic solutions we seek to solve \eqref{e:eq1n=1} and \eqref{e:eq2n=1}.

This whole argument is deceptively simple in the case $n=1$ and hides some very interesting algebra that will be explained in Sections \ref{s:n=2} and \ref{s:3+GeneratorsThm}. 

\begin{rem}
In fact, when $n=1$ none of this is necessary to solve \eqref{e:eq1n=1} and \eqref{e:eq2n=1}, which can be solved via simple linear algebra as above and some observations about holomorphic functions.  Indeed, since we have equations \eqref{e:n=1Identity} and \eqref{e:n=1Identity2}, and impose \eqref{e:NewCCn=1} we can solve for $h_1$ and $k_1$ to get:
\begin{align*}
    h_1 & = \frac{\hat{F}_1}{\left(\hat{F}_1F_1+\hat{G}_1G_1\right)}\\
    k_1 & = -\frac{\hat{G}_1}{\left(\hat{F}_1F_1+\hat{G}_1G_1\right)}.
\end{align*}
By the hypotheses in place it is clear that $h_1, k_1\in H^\infty(\mathbb{D})$ with norms controlled by a function of $\delta$.  Direct algebraic computation shows that for these choices of $h_1$ and $k_1$ that:
\begin{align*}
    F_1h_1-G_1\hat{k}_1 &=1.
\end{align*}
Again, one simply perturbs these solutions in the following way to create smooth and bounded solutions to the problem:
\begin{align}
    H_1(z) & := h_1(z)+\lambda_1(z)G_1(z)\label{e:holo1}\\
    K_1(z) & :=k_1(z)+\hat{\lambda}_1(z)\hat{F}_1(z)\label{e:holo2}.
\end{align}
For any choice of $\lambda_1$ that is smooth and bounded we will have that 
\begin{align*}
    F_1H_1-G_1\hat{K}_1 &=1.
\end{align*}
Observe it is possible to choose the function $\lambda_1$ to force \eqref{e:eq2n=1} to hold too.  Via simple algebra \eqref{e:eq2n=1} is true, namely
$$
F_1K_1+G_1\hat{H}_1=0
$$
if and only if it is possible to solve
$$
F_1k_1+G_1\hat{h}_1+\hat{\lambda}_1\left(\hat{F}_1F_1+\hat{G}_1G_1\right)=0.
$$
By \eqref{e:NewCCn=1} simply choose 
$$
\hat{\lambda}_1:=-\frac{F_1 k_1+G_1\hat{h}_1}{\left(\hat{F}_1F_1+\hat{G}_1G_1\right)}
$$
and this function will be holomorphic and bounded since $F_1$, $G_1$, $h_1$ and $k_1$ are and \eqref{e:NewCCn=1} holds and so division by the function $(\hat{F}_1F_1+\hat{G}_1G_1)$ is allowed.  Then it is immediate that \eqref{e:holo1} and \eqref{e:holo2} belong to $H^\infty(\mathbb{D})$ with norm controlled by a function of $\delta$.
\end{rem}

\subsection{Explanation and Proof in the Case of \texorpdfstring{$n=2$}{n is 2}}
\label{s:n=2}
We next focus on the case of $n=2$ since the resulting algebra is easier to explain and follow, though more difficult than the case of $n=1$.  The modifications to handle the case when $n\geq 3$ (the case of more generators) follows similar, but more involved, arguments.  The reader who is more comfortable with basic algebra connected to Corona problems on the disk can jump to Section \ref{s:3+GeneratorsThm} immediately.

The equations to solve when $n=2$ are:
\begin{align}
\label{e:NewBezout}
    F_1H_1-G_1\hat{K}_1+F_2H_2-G_2\hat{K}_2 &=1\\
\label{e:NewBezout1}
    F_1K_1+G_1\hat{H}_1+F_2K_2+G_2\hat{H}_2 &= 0.
\end{align}
  The goal is to solve \eqref{e:NewBezout} and \eqref{e:NewBezout1} for solutions $H_1, H_2, K_1$, and $K_2\in H^\infty(\mathbb{D})$ with norm control when the functions $F_1, F_2, G_1$ and $G_2$ have no common zeros, and so in particular, we have:
$$
0<\delta^2\leq\vert F_1(z)\vert^2+\vert G_1(z)\vert^2+\vert F_2(z)\vert^2+\vert G_2(z)\vert^2\leq 1 \quad\forall z\in\mathbb{D}.
$$
We split the procedure into smaller steps.

\medskip
{\it Step 1: Key Identities}. 
Observe that for any solutions $h_j$ and $k_j$ the equations \eqref{e:NewBezout} and \eqref{e:NewBezout1} imply the following identities:
\begin{align}
    F_1 & = (\hat{F}_1F_1+\hat{G}_1G_1)\hat{h}_1+(\hat{F}_2F_1+\hat{G}_1G_2)\hat{h}_2+(\hat{G}_1F_2-\hat{G}_2F_1)k_2\label{e:F1identity}\\
    G_1 & = -(\hat{F}_1F_1+\hat{G}_1G_1)k_1+(\hat{F}_2G_1-\hat{F}_1G_2)\hat{h}_2-(\hat{F}_1F_2+\hat{G}_2G_1)k_2\label{e:G1identity}\\
    F_2 & = (\hat{F}_1F_2+\hat{G}_2G_1)\hat{h}_1+(\hat{G}_2F_1-\hat{G}_1F_2)k_1+(\hat{F}_2F_2+\hat{G}_2G_2)\hat{h}_2\label{e:F2identity}\\
    G_2 & = (\hat{F}_1G_2-\hat{F}_2G_1)\hat{h}_1-(\hat{F}_2F_1+\hat{G}_1G_2)k_1-(\hat{F}_2F_2+\hat{G}_2G_2)k_2\label{e:G2identity}.
\end{align}
Indeed, for example to arrive at the identity for $F_1$, \eqref{e:F1identity}, take $\widehat{\cdot}$ of \eqref{e:NewBezout} and then multiply that equation by $F_1$.  Multiply \eqref{e:NewBezout1} by $\hat{G}_1$ and add the resulting equations together.  The identities for $G_1$, $F_2$ and $G_2$ are proved analogously.

One can then substitute these identities into \eqref{e:NewBezout} to arrive at the following auxiliary Bezout equation that is crucial in what follows:
\begin{align}
     (&\hat{F}_1F_1+\hat{G}_1G_1)(\hat{h}_1h_1+\hat{k}_1k_1)+(\hat{F}_2F_1 +\hat{G}_1G_2)(\hat{h}_2h_1+\hat{k}_2k_1)\notag\\
     +(&\hat{G}_1F_2-\hat{G}_2F_1)(h_1k_2-k_1h_2)+(\hat{F}_2G_1-\hat{F}_1G_2)(\hat{h}_1\hat{k}_2-\hat{k}_1\hat{h}_2)\notag\\
      +(&\hat{F}_1F_2+\hat{G}_2G_1)(\hat{h}_1h_2+\hat{k}_1k_2)+(\hat{F}_2F_2+\hat{G}_2G_2)(\hat{h}_2h_2+\hat{k}_2k_2)=1     \label{e:NewBezout2}.
\end{align}

Additionally note that \eqref{e:F1identity}-\eqref{e:G2identity} imply the following by taking $\widehat{\cdot}$ and then $\overline{\partial}$:
\begin{align}
    0 & = (\hat{F}_1F_1+\hat{G}_1G_1)\overline{\partial}h_1+(\hat{F}_1F_2+\hat{G}_2G_1)\overline{\partial} h_2+(\hat{F}_2G_1-\hat{F}_1G_2)\overline{\partial}\hat{k}_2\label{e:dbar1}\\
    0 & = -(\hat{F}_1F_1+\hat{G}_1G_1)\overline{\partial}\hat{k}_1+(\hat{G}_1F_2-\hat{G}_2F_1)\overline{\partial}h_2-(\hat{F}_2F_1+\hat{G}_1G_2)\overline{\partial}\hat{k}_2\label{e:dbar2}\\
    0 & = (\hat{F}_2F_1+\hat{G}_1G_2)\overline{\partial}h_1+(\hat{F}_1G_2-\hat{F}_2G_1)\overline{\partial}\hat{k}_1+(\hat{F}_2F_2+\hat{G}_2G_2)\overline{\partial} h_2\label{e:dbar3}\\
    0 & = (\hat{G}_2F_1-\hat{G}_1F_2)\overline{\partial}h_1-(\hat{F}_1F_2+\hat{G}_2G_1)\overline{\partial}\hat{k}_1-(\hat{F}_2F_2+\hat{G}_2G_2)\overline{\partial}\hat{k}_2\label{e:dbar4}.
\end{align}

\medskip
{\it Step 2: Constructing Smooth and Bounded Solutions.}
We next construct smooth and bounded solutions to each of the equations \eqref{e:NewBezout} and \eqref{e:NewBezout1}.  Namely, we will produce $h_1$, $h_2$, $k_1$ and $k_2$ that are smooth and belong to $L^\infty(\mathbb{D)}$ and additionlly will estimate the norm of the solutions $h_j$ and $k_j$.

First, we solve \eqref{e:NewBezout} as this is the easiest of the equations to satisfy.  This is accomplished by defining
\begin{align*}
    \varphi_1(z) &:=\frac{\overline{F_1(z)}}{D(z)} &    \psi_1(z) &:=-\frac{G_1(\overline{z})}{D(\overline{z})}\\
    \varphi_2(z) &:=\frac{\overline{F_2(z)}}{D(z)} &
    \psi_2(z) &:=-\frac{G_2(\overline{z})}{D(\overline{z})}
\end{align*}
and
$$
D(z)  := \vert F_1(z)\vert^2+\vert G_1(z)\vert^2+\vert F_2(z)\vert^2+\vert G_2(z)\vert^2.
$$
It is clear from the definitions of $\varphi_1,\varphi_2,\psi_1,\psi_2$, that $D(z)\geq\delta^2>0$, and $F_1,F_2,G_1,G_2\in H^\infty(\mathbb{D})$ that we have $\varphi_1,\varphi_2,\psi_1,\psi_2\in L^\infty(\mathbb{D})$ with norm controlled by $\frac{1}{\delta^2}$.  With these functions we immediately have:
\begin{align*}
    F_1(z)\varphi_1(z)-G_1(z)\hat{\psi}_1(z)+F_2(z)\varphi_2(z)-G_2(z)\hat{\psi}_2(z) & =  F_1(z)\frac{\overline{F_1(z)}}{D(z)}-G_1(z)\left(-\frac{\overline{G_1(z)}}{D(z)}\right)\\
    & + F_2(z)\frac{\overline{F_2(z)}}{D(z)}-G_2(z)\left(-\frac{\overline{G_2(z)}}{D(z)}\right)\\
    & = \frac{D(z)}{D(z)}=1
\end{align*}
where we have used that
$$
\hat{\psi}_j(z)=\overline{\psi_j(\overline{z})}=-\frac{\overline{G_j(z)}}{D(z)}\quad \forall j=1,2.
$$
Next, define the following functions:

\begin{align}
    &h_1(z) := \varphi_1(z)+\lambda_1(z)G_2(z)+\lambda_4(z) F_2(z)+\lambda_6(z) G_1(z)\label{e:Smoothfunctions1}
    \\
    &
    k_1(z):= \psi_1(z)+\hat{\lambda}_2(z)\hat{G}_2(z)+\hat{\lambda}_5(z) \hat{F}_2(z)+\hat{\lambda}_6(z)\hat{F}_1(z)\label{e:Smoothfunctions2}
    \\
    &
    h_2(z) := \varphi_2(z)+\lambda_3(z)G_2(z)-\lambda_4(z) F_1(z)+\lambda_5(z) G_1(z)\label{e:Smoothfunctions3}
    \\
    &
    k_2(z) := \psi_2(z)+\hat{\lambda}_1(z)\hat{F}_1(z)-\hat{\lambda}_2(z) \hat{G}_1(z)+\hat{\lambda}_3(z) \hat{F}_2(z).\label{e:Smoothfunctions4}
\end{align}

Here $\lambda_1,\ldots,\lambda_6$ are arbitrary smooth and bounded functions that will be chosen to satisfy certain additional conditions.

It holds that for any choice of $\lambda_1,\ldots,\lambda_6$ we still have solutions to \eqref{e:NewBezout} and \eqref{e:NewBezout1}.  We first observe that with the choice of $h_1, h_2, k_1$ and $k_2$ above, we have a solution to \eqref{e:NewBezout}:
$$
F_1h_1-G_1\hat{k}_1+F_2h_2-G_2\hat{k}_2 = 1.
$$
This is a direct computation.  Indeed, expanding all the algebra gives:
\begin{align*}
    F_1h_1-G_1\hat{k}_1+F_2h_2&-G_2\hat{k}_2
     \\
     &= F_1(\varphi_1+\lambda_1G_2+\lambda_4 F_2+\lambda_6 G_1)-G_1(\hat{\psi}_1+\lambda_2 G_2+\lambda_5 F_2+\lambda_6 F_1)\\
    & + F_2(\varphi_2+\lambda_3G_2-\lambda_4 F_1+\lambda_5 G_1)-G_2(\hat{\psi}_2+\lambda_1 F_1-\lambda_2 G_1+\lambda_3 F_2)\\
    & = (F_1\varphi_1 -G_1\hat{\psi}_1+F_2\varphi_2-G_2\hat{\psi}_2) +F_1\lambda_1G_2+F_1\lambda_4F_2+F_1\lambda_6G_1\\
    & - G_1\lambda_2G_2-G_1\lambda_5F_2-G_1\lambda_6F_1+F_2\lambda_3G_2-F_2\lambda_4 F_1+F_2\lambda_5G_1\\
    &-G_2\lambda_1F_1+G_2\lambda_2G_1-G_2\lambda_3F_2\\
    & = 1.
\end{align*}

Next, impose an additional condition on the $\lambda_1,\ldots,\lambda_6$ to force the solutions to satisfy \eqref{e:NewBezout1}:
$$
F_1k_1+G_1\hat{h}_1+F_2k_2+G_2\hat{h}_2 = 0.
$$
Using the ansatz for the solutions $h_1,h_2,k_1$ and $k_2$ we have
\begin{align*}
    F_1k_1&+G_1\hat{h}_1+F_2k_2+G_2\hat{h}_2
     \\
     &= F_1(\psi_1+\hat{\lambda}_2\hat{G}_2+\hat{\lambda}_5 \hat{F}_2+\hat{\lambda}_6\hat{F}_1)+G_1(\hat{\varphi}_1+\hat{\lambda}_1\hat{G}_2+\hat{\lambda}_4 \hat{F}_2+\hat{\lambda}_6 \hat{G}_1)\\
    & + F_2(\psi_2+\hat{\lambda}_1\hat{F}_1-\hat{\lambda}_2 \hat{G}_1+\hat{\lambda}_3 \hat{F}_2)+G_2(\hat{\varphi}_2+\hat{\lambda}_3\hat{G}_2-\hat{\lambda}_4 \hat{F}_1+\hat{\lambda}_5 \hat{G}_1)\\
    & = (F_1\psi_1+G_1\hat{\varphi}_1+F_2\psi_2+G_2\hat{\varphi}_2) +\hat{\lambda}_1(G_1\hat{G}_2+\hat{F}_1F_2)+\hat{\lambda}_2(F_1\hat{G}_2-F_2\hat{G}_1)\\
    & +\hat{\lambda}_3(F_2\hat{F}_2+G_2\hat{G}_2)+\hat{\lambda}_4(G_1\hat{F_2}-\hat{F}_1G_2)+\hat{\lambda}_5(F_1\hat{F}_2+\hat{G}_1G_2)+\hat{\lambda}_6(F_1\hat{F}_1+G_1\hat{G}_1).
\end{align*}
But, from the definitions of $\psi_1,\psi_2,\varphi_1,\varphi_2$ we have:
$$
F_1(z)\psi_1(z)+G_1(z)\hat{\varphi}_1(z)+F_2(z)\psi_2(z)+G_2(z)\hat{\varphi}_2(z)=\frac{1}{D(\overline{z})}\sum_{j=1}^{2}\left( F_j(\overline{z})G_j(z)-F_j(z)G_j(\overline{z})\right).
$$
Focus on the right hand side and define:
$$
\tau(z):=\frac{1}{D(\overline{z})}\sum_{j=1}^{2}\left( F_j(\overline{z})G_j(z)-F_j(z)G_j(\overline{z})\right)
$$
and observe that $\tau\in L^\infty(\mathbb{D})$ because by hypothesis $D(z)\geq\delta^2>0$ and $F_j, G_j\in H^\infty(\mathbb{D})$ and further $\Vert \tau\Vert_{L^\infty(\mathbb{D})}\lesssim \frac{1}{\delta^2}$.
Hence, we have that \eqref{e:NewBezout1} holds, i.e.,
$$
F_1k_1+G_1\hat{h}_1+F_2k_2+G_2\hat{h}_2 = 0
$$
if and only if it is possible to solve:
\begin{align}
    \label{e:lambdaconstraints}
    -\tau & = \hat{\lambda}_1(G_1\hat{G}_2+\hat{F}_1F_2)+\hat{\lambda}_2(F_1\hat{G}_2-F_2\hat{G}_1)+\hat{\lambda}_3(F_2\hat{F}_2+G_2\hat{G}_2)\\
    & +\hat{\lambda}_4(G_1\hat{F_2}-\hat{F}_1G_2)+\hat{\lambda}_5(F_1\hat{F}_2+\hat{G}_1G_2)+\hat{\lambda}_6(F_1\hat{F}_1+G_1\hat{G}_1)\notag
\end{align}
for the functions $\lambda_1,\ldots,\lambda_6$.

Define
\begin{align}
\label{e:NewDenominator}
\Delta(z) & :=\abs{G_1(z)\hat{G}_2(z)+\hat{F}_1(z)F_2(z)}^2+\abs{F_1(z)\hat{G}_2(z)-F_2(z)\hat{G}_1(z)}^2
\\
&
+
\abs{F_2(z)\hat{F}_2(z)+G_2(z)\hat{G}_2(z)}^2
+\abs{G_1(z)\hat{F_2}(z)-\hat{F}_1(z)G_2(z)}^2
\notag\\
&
+\abs{F_1(z)\hat{F}_2(z)+\hat{G}_1(z)G_2(z)}^2+\abs{F_1(z)\hat{F}_1(z)+G_1(z)\hat{G}_1(z)}^2.\notag
\end{align}
It is the case that $\Delta(z)\geq\delta^2>0$; this follows immediately from \eqref{e:NewBezout2} and is a quantitative version of an observation of Gentili, Sarfatti and Struppa \cite{MR3621101}*{Theorem 3.1, p. 1654}.  In particular we remark at this point that the argument given shows that a condition like $\Delta(z)\geq \delta^2>0$ is necessary for the possible solutions to the the Theorem \ref{t:NewCoronaS2}.

Hence, it is then possible to solve \eqref{e:lambdaconstraints} for $\lambda_1,\ldots,\lambda_6$ with
\begin{align*}
    \hat{\lambda}_1(z) & :=-\frac{\overline{G_1(z)\hat{G}_2(z)+\hat{F}_1(z)F_2(z)}}{\Delta(z)}\tau(z) & \quad
    \hat{\lambda}_2(z) & :=-\frac{\overline{F_1(z)\hat{G}_2(z)-F_2(z)\hat{G}_1(z)}}{\Delta(z)}\tau(z)\\
    \hat{\lambda}_3(z) & :=-\frac{\overline{F_2(z)\hat{F}_2(z)+G_2(z)\hat{G}_2(z)}}{\Delta(z)}\tau(z) & \quad
    \hat{\lambda}_4(z) & :=-\frac{\overline{G_1(z)\hat{F_2}(z)-\hat{F}_1(z)G_2(z)}}{\Delta(z)}\tau(z)\\
    \hat{\lambda}_5(z) & :=-\frac{\overline{F_1(z)\hat{F}_2(z)+\hat{G}_1(z)G_2(z)}}{\Delta(z)}\tau(z) & \quad
    \hat{\lambda}_6(z) & :=-\frac{\overline{F_1(z)\hat{F}_1(z)+G_1(z)\hat{G}_1(z)}}{\Delta(z)}\tau(z).
\end{align*}
Note that we can solve for $\lambda_j$ by using that $p(z)=\hat{\hat{p}}(z)$.  It is clear from the definitions of $\lambda_1,\ldots,\lambda_6$ that each of these functions are in $L^\infty(\mathbb{D})$ with norm controlled by $\frac{1}{\delta^4}$.  Since we can solve \eqref{e:lambdaconstraints}, we have solved \eqref{e:NewBezout1}.  Then with this choice of $\lambda_1,\ldots,\lambda_6$ we have solved \eqref{e:NewBezout} and \eqref{e:NewBezout1} with bounded and smooth solutions $h_1$, $h_2$, $k_1$ and $k_2$ with $\Vert h_j\Vert_{L^\infty(\mathbb{D})}~\lesssim ~\frac{1}{\delta^2}+~\frac{1}{\delta^4}$ and $\Vert k_j\Vert_{L^\infty(\mathbb{D})}\lesssim \frac{1}{\delta^2}+\frac{1}{\delta^4}$.

\medskip
{\it Step 3: Perturbing the Solutions to be Holomorphic.}
We next perturb the smooth and bounded solutions (\ref{e:Smoothfunctions1})-(\ref{e:Smoothfunctions4}) to solve (\ref{e:NewBezout}) and (\ref{e:NewBezout1}) in a holomorphic manner.  
\\
One observes that
\begin{align*}
H_1(z) &:= h_1(z)+\beta_1(z)(\hat{F}_1(z)F_2(z)+\hat{G}_2(z)G_1(z))+\beta_2(z)(\hat{F}_2(z)F_2(z)+\hat{G}_2(z)G_2(z))\\
&- \beta_3(z)(\hat{F}_2(z)G_1(z)-\hat{F}_1(z)G_2(z))\\
K_1(z) & := k_1(z)+\hat{\beta}_1(z)(\hat{F}_1(z)G_2(z)-\hat{F}_2(z)G_1(z)) - \hat{\beta}_3(z)(\hat{F}_1(z)F_2(z)+\hat{G}_2(z)G_1(z))\\
& -\hat{\beta}_4(z)(\hat{F}_2(z)F_2(z)+\hat{G}_2(z)G_2(z))\\
H_2(z) & := h_2(z)-\beta_1(z)(\hat{F}_1(z)F_1(z)+\hat{G}_1(z)G_1(z))-\beta_2(z)(\hat{F}_2(z)F_1(z)+\hat{G}_1(z)G_2(z))\\
& +\beta_4(z)(\hat{F}_1(z)G_2(z)-\hat{F}_2(z)G_1(z))\\
K_2(z) & := k_2(z)+\hat{\beta}_2(z)(\hat{F}_1(z)G_2(z)-\hat{F}_2(z)G_1(z)) + \hat{\beta}_3(z)(\hat{F}_1(z)F_1(z)+\hat{G}_1(z)G_1(z))\\
& + \hat{\beta}_4(z)(\hat{F}_2(z)F_1(z)+\hat{G}_1(z)G_2(z))
\end{align*}
are solutions to \eqref{e:NewBezout} and \eqref{e:NewBezout1} for any choice of $\beta_1,\ldots, \beta_4$.  Indeed, we have that for any $\beta_1,\ldots,\beta_4$ that $H_1, H_2, K_1$ and $K_2$ are solutions to the system:
\begin{align*}
    F_1H_1-G_1\hat{K}_1+F_2H_2-G_2\hat{K}_2 &=1\\
    F_1K_1+G_1\hat{H}_1+F_2K_2+G_2\hat{H}_2 &= 0,
\end{align*}
because the $\beta_1,\ldots, \beta_4$ arise a solutions to the related homogeneous system.  We verify that these are solutions now via direct algebra.  For the first equation we have:
\begin{align*}
F_1H_1-G_1\hat{K}_1&+F_2H_2-G_2\hat{K}_2
\\
&
 = F_1\left(h_1+\beta_1(\hat{F}_1F_2+\hat{G}_2G_1)+\beta_2(\hat{F}_2F_2+\hat{G}_2G_2)- \beta_3(\hat{F}_2G_1-\hat{F}_1G_2)\right)\\
& -G_1\left(\hat{k}_1+\beta_1(F_1\hat{G}_2-F_2\hat{G}_1)-\beta_3(F_1\hat{F}_2+G_2\hat{G}_1) -\beta_4(F_2\hat{F}_2+G_2\hat{G}_2)\right)\\
& + F_2\left(h_2-\beta_1(\hat{F}_1F_1+\hat{G}_1G_1)-\beta_2(\hat{F}_2F_1+\hat{G}_1G_2)-\beta_4(\hat{F}_2G_1-\hat{F}_1G_2)\right)\\
& -G_2\left(\hat{k}_2+\beta_2(F_1\hat{G}_2-F_2\hat{G}_1) + \beta_3(F_1\hat{F}_1+G_1\hat{G}_1) + \beta_4(F_2\hat{F}_1+G_1\hat{G}_2)\right)\\
& = F_1h_1-G_1\hat{k}_1+F_2h_2-G_2\hat{k}_2\\
& + \beta_1\left(F_1(\hat{F}_1F_2+\hat{G}_2G_1)-G_1(F_1\hat{G}_2-F_2\hat{G}_1)-F_2(\hat{F}_1F_1+\hat{G}_1G_1)\right)\\
& + \beta_2\left(F_1(\hat{F}_2F_2+\hat{G}_2G_2)-F_2(\hat{F}_2F_1+\hat{G}_1G_2)-G_2(F_1\hat{G}_2-F_2\hat{G}_1)\right)\\
& + \beta_3\left(-F_1(\hat{F}_2G_1-\hat{F}_1G_2)+G_1(F_1\hat{F}_2+G_2\hat{G}_1)-G_2(F_1\hat{F}_1+G_1\hat{G}_1)\right)\\
& + \beta_4\left(G_1(F_2\hat{F}_2+G_2\hat{G}_2)-F_2(\hat{F}_2G_1-\hat{F}_1G_2)-G_2(F_2\hat{F}_1+G_1\hat{G}_2)\right)\\
& = 1.
\end{align*}
The last equality holds by the properties of the solutions $h_1,h_2,k_1$ and $k_2$ and one can check by inspection that each of the factors multiplying a $\beta_j$ sum to zero.

Similarly, checking the second equation yields:
\begin{align*}
    F_1K_1+G_1\hat{H}_1&+F_2K_2+G_2\hat{H}_2
    \\
    &
    = F_1\left(k_1+\hat{\beta}_1(\hat{F}_1G_2-\hat{F}_2G_1) - \hat{\beta}_3(\hat{F}_1F_2+\hat{G}_2G_1)-\hat{\beta}_4(\hat{F}_2F_2+\hat{G}_2G_2)\right)\\
    & + G_1\left(\hat{h}_1+\hat{\beta}_1(F_1\hat{F}_2+G_2\hat{G}_1)+\hat{\beta}_2(F_2\hat{F}_2+G_2\hat{G}_2) - \hat{\beta}_3(F_2\hat{G}_1-F_1\hat{G}_2)\right)\\
    & + F_2\left(k_2+\hat{\beta}_2(\hat{F}_1G_2-\hat{F}_2G_1) + \hat{\beta}_3(\hat{F}_1F_1+\hat{G}_1G_1) + \hat{\beta}_4(\hat{F}_2F_1+\hat{G}_1G_2)\right)\\
    & + G_2\left(\hat{h}_2-\hat{\beta}_1(F_1\hat{F}_1+G_1\hat{G}_1)-\hat{\beta}_2(F_2\hat{F}_1+G_1\hat{G}_2)-\hat{\beta}_4(F_2\hat{G}_1-F_1\hat{G}_2)\right)\\
    & = F_1k_1+G_1\hat{h}_1+F_2k_2+ G_2\hat{h}_2\\
    & +\hat{\beta}_1\left(F_1(\hat{F}_1G_2-\hat{F}_2G_1)+G_1(F_1\hat{F}_2+G_2\hat{G}_1)-G_2(F_1\hat{F}_1+G_1\hat{G}_1)\right)\\
    & +\hat{\beta}_2\left(G_1(F_2\hat{F}_2+G_2\hat{G}_2)+F_2(\hat{F}_1G_2-\hat{F}_2G_1)-G_2(F_2\hat{F}_1+G_1\hat{G}_2)\right)\\
    & +\hat{\beta}_3\left(-F_1(\hat{F}_1F_2+\hat{G}_2G_1)-G_1(F_2\hat{G}_1-F_1\hat{G}_2)+F_2(\hat{F}_1F_1+\hat{G}_1G_1)\right)\\
    & +\hat{\beta}_4\left(-F_1(\hat{F}_2F_2+\hat{G}_2G_2)+F_2(\hat{F}_2F_1+\hat{G}_1G_2)-G_2(F_2\hat{G}_1-F_1\hat{G}_2)\right)\\
    & =0.
\end{align*}
Again, last equality holds by the properties of the solutions $h_1,h_2,k_1$ and $k_2$ and one can check by inspection that each of the factors multiplying a $\beta_j$ again sum to zero.

Finally, we choose $\beta_1,\ldots,\beta_4$ so that the functions $H_1,K_1,H_2,K_2$ are holomorphic.  Observe that for appropriate choice of $\beta_j$ we have
\begin{align*}
    0=\overline{\partial}H_1 & = \overline{\partial}h_1+(\hat{F}_1F_2+\hat{G}_2G_1)\overline{\partial}\beta_1+(\hat{F}_2F_2+\hat{G}_2G_2)\overline{\partial}\beta_2- (\hat{F}_2G_1-\hat{F}_1G_2)\overline{\partial}\beta_3\\
    0=\overline{\partial}K_1 & =\overline{\partial}k_1+ (\hat{F}_1G_2-\hat{F}_2G_1)\overline{\partial}\hat{\beta}_1 - (\hat{F}_1F_2+\hat{G}_2G_1)\overline{\partial}\hat{\beta}_3-(\hat{F}_2F_2+\hat{G}_2G_2)\overline{\partial}\hat{\beta}_4\\
    0=\overline{\partial}H_2 & = \overline{\partial}h_2-(\hat{F}_1F_1+\hat{G}_1G_1)\overline{\partial}\beta_1-(\hat{F}_2F_1+\hat{G}_1G_2)\overline{\partial}\beta_2 + (\hat{F}_1G_2-\hat{F}_2G_1)\overline{\partial}\beta_4\\
    0=\overline{\partial}K_2 & = \overline{\partial}k_2+(\hat{F}_1G_2-\hat{F}_2G_1)\overline{\partial}\hat{\beta}_2 + (\hat{F}_1F_1+\hat{G}_1G_1)\overline{\partial}\hat{\beta}_3 + (\hat{F}_2F_1+\hat{G}_1G_2)\overline{\partial}\hat{\beta}_4.
\end{align*}

In particular, choosing the functions $\beta_j$ to satisfy:
\begin{align*}
    \overline{\partial}\beta_1 & = \overline{\partial}h_2 (\hat{h}_1h_1+\hat{k}_1k_1) -\overline{\partial}h_1 (\hat{h}_1h_2 +\hat{k}_1k_2) - \overline{\partial}\hat{k}_1 (h_2k_1-h_1k_2)
    \\
    \overline{\partial}\beta_2 & = \overline{\partial}h_2 (\hat{h}_2h_1 + \hat{k}_2k_1) - \overline{\partial}h_1(\hat{h}_2h_2+\hat{k}_2k_2) + \overline{\partial}\hat{k}_2(h_1k_2-h_2k_1)
    \\
     {\overline{\partial}\beta_3} & =  { -\overline{\partial}h_1(\hat{h}_2\hat{k}_1-\hat{h}_1\hat{k}_2)- \overline{\partial}\hat{k}_2(\hat{h}_1h_1+\hat{k}_1k_1) + \overline{\partial}\hat{k}_1 (\hat{h}_2h_1+\hat{k}_2k_1)}
    \\
    \overline{\partial}\beta_4 & = \overline{\partial}h_2(\hat{h}_1\hat{k}_2-\hat{h}_2\hat{k}_1)- \overline{\partial}\hat{k}_2(\hat{h}_1h_2+\hat{k}_1k_2) + \overline{\partial}\hat{k}_1(\hat{h}_2h_2+\hat{k}_2k_2).
\end{align*}

With these choice of $\beta_j$ it is possible to have $H_1,K_1,H_2$ and $K_2$ being holomorphic.  We include only the computations to verify this for the term $H_1$ and $K_2$ as the other cases are handled identically.  The idea is to substitute in the formulas for $\overline{\partial}\beta_j$, collect terms relative to the expression we are trying to cancel, use identities \eqref{e:dbar1}, \eqref{e:dbar2}, \eqref{e:dbar3}, \eqref{e:dbar4} to convert some terms to the ones needed, and then use \eqref{e:NewBezout} to show that things cancel.

Turning first to demonstrate the holomorphicity of $H_1$, we have
\begin{align*}
    \overline{\partial}H_1 & = \overline{\partial}h_1+(\hat{F}_1F_2+\hat{G}_2G_1)\overline{\partial}\beta_1
    +(\hat{F}_2F_2+\hat{G}_2G_2)\overline{\partial}\beta_2- (\hat{F}_2G_1-\hat{F}_1G_2)\overline{\partial}\beta_3
    \\
    &
     = \overline{\partial}h_1+(\hat{F}_1F_2+\hat{G}_2G_1)\left(\overline{\partial}h_2 (\hat{h}_1h_1+\hat{k}_1k_1) -\overline{\partial}h_1 (\hat{h}_1h_2 +\hat{k}_1k_2) - \overline{\partial}\hat{k}_1 (h_2k_1-h_1k_2)\right)
     \\
    &
    + (\hat{F}_2F_2+\hat{G}_2G_2)\left(\overline{\partial}h_2 (\hat{h}_2h_1 + \hat{k}_2k_1) - \overline{\partial}h_1(\hat{h}_2h_2+\hat{k}_2k_2) + \overline{\partial}\hat{k}_2(h_1k_2-h_2k_1)\right)
    \\
    &
     -(\hat{F}_2G_1-\hat{F}_1G_2)
      {\left(    -\overline{\partial}h_1(\hat{h}_2\hat{k}_1-\hat{h}_1\hat{k}_2)- \overline{\partial}\hat{k}_2(\hat{h}_1h_1+\hat{k}_1k_1) + \overline{\partial}\hat{k}_1 (\hat{h}_2h_1+\hat{k}_2k_1)\right)}
     \\
    &
    = \overline{\partial}h_1\left(1-(\hat{F}_1F_2+\hat{G}_2G_1)(\hat{h}_1h_2 +\hat{k}_1k_2)-(\hat{F}_2F_2+\hat{G}_2G_2)(\hat{h}_2h_2+\hat{k}_2k_2)\right.
    \\
    &
     \left.-(\hat{F}_2G_1-\hat{F}_1G_2)(\hat{h}_1\hat{k}_2-\hat{h}_2\hat{k}_1)
    \right)
    \\
    &
     +(\hat{h}_1h_1+\hat{k}_1k_1)\left((\hat{F}_1F_2+\hat{G}_2G_1)\overline{\partial}h_2+(\hat{F}_2G_1-\hat{F}_1G_2)\overline{\partial}\hat{k}_2\right)
     \\
    &
    + (h_2k_1-h_1k_2)\left(-(\hat{F}_1F_2+\hat{G}_2G_1)\overline{\partial}\hat{k}_1-(\hat{F}_2F_2+\hat{G}_2G_2)\overline{\partial}\hat{k}_2\right)
    \\
    &
    + (\hat{h}_2h_1 + \hat{k}_2k_1)\left((\hat{F}_2F_2+\hat{G}_2G_2)\overline{\partial}h_2-(\hat{F}_2G_1-\hat{F}_1G_2)\overline{\partial}\hat{k}_1 \right)
    \end{align*}
    and continuing,
    \begin{align*}
    \overline{\partial}H_1 & =
     \overline{\partial}h_1\left(1-(\hat{F}_1F_2+\hat{G}_2G_1)(\hat{h}_1h_2 +\hat{k}_1k_2)-(\hat{F}_2F_2+\hat{G}_2G_2)(\hat{h}_2h_2+\hat{k}_2k_2)\right.
    \\
     & \left.-(\hat{F}_2G_1-\hat{F}_1G_2)(\hat{h}_1\hat{k}_2-\hat{h}_2\hat{k}_1)\right)
     -(\hat{F}_1F_1+\hat{G}_1G_1)(\hat{h}_1h_1+\hat{k}_1k_1)\overline{\partial}h_1
     \\
     &
     -(\hat{G}_2F_1-\hat{G}_1F_2)(h_2k_1-h_1k_2)\overline{\partial}h_1
     -(\hat{F}_2F_1+\hat{G}_1G_2)(\hat{h}_2\hat{h}_1
     +\hat{k}_2k_1)\overline{\partial}h_1
     \\
     &
     = \overline{\partial}h_1\Big(1-(\hat{F}_1F_2+\hat{G}_2G_1)(\hat{h}_1h_2 +\hat{k}_1k_2)-(\hat{F}_2F_2+\hat{G}_2G_2)(\hat{h}_2h_2+\hat{k}_2k_2)
     \\
     &
     -(\hat{F}_2G_1-\hat{F}_1G_2)(\hat{h}_1\hat{k}_2-\hat{h}_2\hat{k}_1)
     -(\hat{F}_1F_1+\hat{G}_1G_1)(\hat{h}_1h_1+\hat{k}_1k_1)
     \\
     &
     -(\hat{F}_2F_1+\hat{G}_1G_2)(\hat{h}_2h_1+\hat{k}_2k_1)
     -(\hat{G}_2F_1-\hat{G}_1F_2)(h_2k_1-h_1k_2)\Big)
     \\
     & = 0.
\end{align*}
The justification of the algebra is explained in the following way.  The first equality is definition; the second is substituting in the expressions for the $\overline{\partial}\beta_j$.  The next equality is regrouping terms, either collecting $\overline{\partial}h_1$ terms or terms with common factors but no $\overline{\partial}h_1$.  The next equality is using \eqref{e:dbar1} for the term involving $(\hat{h}_1h_1+\hat{k}_1k_1)$, \eqref{e:dbar3} for the term involving $(\hat{h}_2h_1+\hat{k}_2k_1)$, and \eqref{e:dbar4} for the term involving $(h_2k_1-h_1k_2)$. Regrouping and using \eqref{e:NewBezout2} gives the last line.

We next verify the holomorphicity of the function $K_2$.  We have
\begin{align*}
    \overline{\partial}\hat{K}_2 & = \overline{\partial}\hat{k}_2+(\hat{G}_2F_1-\hat{G}_1F_2)\overline{\partial}\beta_2 + (\hat{F}_1F_1+\hat{G}_1G_1)\overline{\partial}\beta_3 + (\hat{F}_1F_2+\hat{G}_2G_1)\overline{\partial}\beta_4\\
    & = \overline{\partial}\hat{k}_2+(\hat{G}_2F_1-\hat{G}_1F_2)\left(    \overline{\partial}h_2 (\hat{h}_2h_1 + \hat{k}_2k_1) - \overline{\partial}h_1(\hat{h}_2h_2+\hat{k}_2k_2) + \overline{\partial}\hat{k}_2(h_1k_2-h_2k_1)\right)\\
    & + (\hat{F}_1F_1+\hat{G}_1G_1)\left(    -\overline{\partial}h_1(\hat{h}_2\hat{k}_1-\hat{h}_1\hat{k}_2)- \overline{\partial}\hat{k}_2(\hat{h}_1h_1+\hat{k}_1k_1) + \overline{\partial}\hat{k}_1 (\hat{h}_2h_1+\hat{k}_2k_1)\right)\\
    & +(\hat{F}_1F_2+\hat{G}_2G_1)\left(\overline{\partial}h_2(\hat{h}_1\hat{k}_2-\hat{h}_2\hat{k}_1)- \overline{\partial}\hat{k}_2(\hat{h}_1h_2+\hat{k}_1k_2) + \overline{\partial}\hat{k}_1(\hat{h}_2h_2+\hat{k}_2k_2)
\right)
\end{align*}
and also
\begin{align*}
    \overline{\partial}\hat{K}_2 & =
\overline{\partial}\hat{k}_2\left(1-(\hat{G}_1F_2-\hat{G}_2F_1)(h_1k_2 -h_2k_1)-(\hat{F}_1F_1+\hat{G}_1G_1)(\hat{h}_1h_1+\hat{k}_1k_1)\right.\\
    & \left.-(\hat{F}_1F_2+\hat{G}_2G_1)(\hat{h}_1h_2+\hat{k}_1k_2)\right)
    \\
    &
     +(\hat{h}_2h_1+\hat{k}_2k_1)\left((\hat{G}_2F_1-\hat{G}_1F_2)
     \overline{\partial}h_2+(\hat{F}_1F_1+\hat{G}_1G_1)\overline{\partial}\hat{k}_1\right)\\
    & + (\hat{h}_2h_2+\hat{k}_2k_2)\left(-(\hat{G}_2F_1-\hat{G}_1F_2)\overline{\partial}h_1+(\hat{F}_1F_2+\hat{G}_2G_1)\overline{\partial}\hat{k}_1\right)\\
    & + (\hat{h}_2\hat{k}_1- \hat{h}_1\hat{k}_2)\left(-(\hat{F}_1F_1+\hat{G}_1G_1)\overline{\partial}h_1-(\hat{F}_1F_2+\hat{G}_2G_1)\overline{\partial}h_2 \right)\\
    & = \overline{\partial}\hat{k}_2\left(1-(\hat{G}_1F_2-\hat{G}_2F_1)(h_1k_2 -h_2k_1)-(\hat{F}_1F_1+\hat{G}_1G_1)(\hat{h}_1h_1+\hat{k}_1k_1)\right.\\
    & \left.-(\hat{F}_1F_2+\hat{G}_2G_1)(\hat{h}_1h_2+\hat{k}_1k_2)\right)-(\hat{F}_2F_1+\hat{G}_1G_2)(\hat{h}_2h_1+\hat{k}_2k_1)\overline{\partial}\hat{k}_2\\
     &  (\hat{F}_2F_2+\hat{G}_2G_2)(\hat{h}_2h_2+\hat{k}_2k_2)
     \overline{\partial}\hat{k}_2-(\hat{F}_1G_2-\hat{F}_2G_1)
     (\hat{h}_2\hat{k}_1-\hat{h}_1\hat{k}_2)\overline{\partial}\hat{k}_2
     \end{align*}
     finally we get
     \begin{align*}
    \overline{\partial}\hat{K}_2 & =  \overline{\partial}\hat{k}_2\left(1-(\hat{G}_1F_2-\hat{G}_2F_1)(h_1k_2 -h_2k_1)-(\hat{F}_1F_1+\hat{G}_1G_1)(\hat{h}_1h_1+\hat{k}_1k_1)\right.
     \\
     &  {-}\left.(\hat{F}_2F_1+\hat{G}_1G_2)(\hat{h}_2h_1+\hat{k}_2k_1)
     -(\hat{F}_2F_2+\hat{G}_2G_2)(\hat{h}_2h_2+\hat{k}_2k_2)\right.
     \\
     &
     -\left.( {\hat{F}_1G_2-\hat{F}_2G_1)(\hat{h}_2\hat{k}_1
     -\hat{h}_1\hat{k}_2)}\right)\\
     & = 0.
\end{align*}
Again, this is justified in the following manner. The first equality is definition; the second is substituting in the expressions for the $\overline{\partial}\beta_j$.  The next equality is regrouping terms, either collecting $\overline{\partial}\hat{k}_2$ terms or terms with common factors but no $\overline{\partial}\hat{k}_2$.  The next equality is using \eqref{e:dbar2} for the term involving $(\hat{h}_2h_1+\hat{k}_2k_1)$, \eqref{e:dbar4} for the term involving $(\hat{h}_2h_2+\hat{k}_2k_2)$, and \eqref{e:dbar1} for the term involving $(\hat{h}_2\hat{k}_1-\hat{h}_1\hat{k}_2)$. Regrouping and using \eqref{e:NewBezout2} gives the last line.

Using the Corona Theorem of Carleson, or more precisely Wolff's proof of estimates to $\overline{\partial}$-equations on the disk, see for example \cite{MR2261424}*{Chapter 8, Theorem 1.2}, we have that the functions $H_1,H_2,K_1,K_2$ are in $H^\infty(\mathbb{D})$ with bound controlled by $C(\delta)$ and they satisfy the equation \eqref{e:NewBezout} and \eqref{e:NewBezout1}.  This will be explained when handling the case of $n$ generators.  See Theorem \ref{t:boundedsolutions} in Section \ref{s:Wolff} for the justification as to why there are bounded holomorphic solutions. This ends the proof in the case $n=2$. 

\subsection{The Case of More Generators: The Proof of Theorem \ref{t:NewCorona}}
\label{s:3+GeneratorsThm}

We turn to the proof of the main theorem. 

We seek to solve the following system of equations:
\begin{align}
    \sum_{j=1}^n F_jH_j-\sum_{j=1}^nG_j\hat{K}_j & =1\label{e:nNewBezout1}\\
    \sum_{j=1}^n F_jK_j+\sum_{j=1}^nG_j\hat{H}_j & =0\label{e:nNewBezout2}.
\end{align}
We will work under the hypothesis that $F_j$ and $G_j$ have no common zeros in a quantitative sense, and hence satisfy the Carleson condition:
\begin{equation}
    \label{e:nCarl}
    0<\delta^2\leq\sum_{j=1}^n \abs{F_j(z)}^2+\sum_{j=1}^n \abs{G_j(z)}^2\leq 1\quad\forall z\in\mathbb{D}.
\end{equation}
We will see in the course of this argument that \eqref{e:NewCC2} is necessary as well.  Again this will be accomplished in three steps mirroring the case of $n=2$, but with suitable modifications to accommodate the additional algebraic challenges.

\subsubsection{Step 1: Key Identities}

Equations \eqref{e:nNewBezout1} and \eqref{e:nNewBezout2} imply a third Bezout equation holds, which is:
\begin{align}
    \label{e:nNewBezout3}
1 &= \sum_{s=1}^{n}\sum_{r=1}^n \left(\hat{F}_rF_s+\hat{G}_sG_r\right)\left(\hat{H}_rH_s+\hat{K}_rK_s\right)+\sum_{s=1}^{n}\sum_{r=s+1}^n \left(\hat{G}_sF_r-\hat{G}_rF_s\right)\left(K_rH_s-K_sH_r\right)\\
    & +\sum_{s=1}^{n}\sum_{r=s+1}^n \left(\hat{F}_rG_s-\hat{F}_sG_r\right)\left(\hat{H}_s\hat{K}_r-\hat{H}_r\hat{K}_s\right).\notag
\end{align}

To see this, take $\widehat{\cdot}$ of \eqref{e:nNewBezout1} and rearrange to see that we have the two equations, with the second equation being \eqref{e:nNewBezout2}:
\begin{align*}
    \sum_{j=1}^n \hat{F}_j\hat{H}_j-\sum_{j=1}^n\hat{G}_jK_j & =1\\
    \sum_{j=1}^nG_j\hat{H}_j+\sum_{j=1}^n F_jK_j & =0.
\end{align*}
Multiply the first question by $F_{s}$ and the second equation by $\hat{G}_{s}$ and add them together to get the identity:
\begin{equation}
\label{e:Fjidentity}
F_{s}=\sum_{r=1}^n \left(\hat{F}_rF_{s}+\hat{G}_{s}G_r\right)\hat{H}_r + \sum_{r=1}^n\left(\hat{G}_{s}F_r-\hat{G}_rF_{s}\right)K_r\quad \forall 1\leq s\leq n;
\end{equation}
and multiply the first equation by $G_{s}$ and the second equation by ${-}\hat{F}_{s}$ and then add the resulting equations to get:
\begin{equation}
\label{e:Gjidentity}
G_{s}=\sum_{r=1}^n \left(\hat{F}_rG_{s}-\hat{F}_{s}G_r\right)\hat{H}_r - \sum_{r=1}^n\left(\hat{F}_{s}F_r+\hat{G}_rG_{s}\right)K_r\quad \forall 1\leq s\leq n.
\end{equation}
Substitute \eqref{e:Fjidentity} and \eqref{e:Gjidentity} into \eqref{e:nNewBezout1} to see that we have:

\begin{align*}
    1 & = \sum_{s=1}^n F_sH_s-\sum_{s=1}^nG_s\hat{K}_s\\
    & = \sum_{s=1}^n\left(\sum_{r=1}^n \left(\hat{F}_rF_{s}+\hat{G}_{s}G_r\right)\hat{H}_r + \sum_{r=1}^n\left(\hat{G}_{s}F_r-\hat{G}_rF_{s}\right)K_r\right)H_s\\
    & -\sum_{s=1}^n\left(\sum_{r=1}^n \left(\hat{F}_rG_{s}-\hat{F}_{s}G_r\right)\hat{H}_r - \sum_{r=1}^n\left(\hat{F}_{s}F_r+\hat{G}_rG_{s}\right)K_r\right)\hat{K}_s\\
    & = \sum_{s=1}^n\sum_{r=1}^n\left(\hat{F}_rF_{s}+\hat{G}_{s}G_r\right)\hat{H}_rH_s+\sum_{s=1}^n\sum_{r=1}^n\left(\hat{F}_sF_{r}+\hat{G}_{r}G_s\right)\hat{K}_sK_r\\
    & +\sum_{s=1}^n\sum_{r=1}^n\left(\hat{G}_{s}F_r-\hat{G}_rF_{s}\right)K_rH_s
    -\sum_{s=1}^n\sum_{r=1}^n\left(\hat{F}_{r}G_s-\hat{F}_sG_{r}\right)\hat{H}_r\hat{K}_s
    \end{align*}
\begin{align*}
& = \sum_{s=1}^{n}\sum_{r=1}^n \left(\hat{F}_rF_s+\hat{G}_sG_r\right)\left(\hat{H}_rH_s+\hat{K}_rK_s\right)+\frac{1}{2}\sum_{s=1}^{n}\sum_{r=1}^n \left(\hat{G}_sF_r-\hat{G}_rF_s\right)\left(K_rH_s-K_sH_r\right)\\
    & +\frac{1}{2}\sum_{s=1}^{n}\sum_{r=1}^n \left(\hat{F}_rG_s-\hat{F}_sG_r\right)\left(\hat{H}_s\hat{K}_r-\hat{H}_r\hat{K}_s\right)\\
    & = \sum_{s=1}^{n}\sum_{r=1}^n \left(\hat{F}_rF_s+\hat{G}_sG_r\right)\left(\hat{H}_rH_s+\hat{K}_rK_s\right)+\sum_{s=1}^{n}\sum_{r=s+1}^n \left(\hat{G}_sF_r-\hat{G}_rF_s\right)\left(K_rH_s-K_sH_r\right)\\
    & +\sum_{s=1}^{n}\sum_{r=s+1}^n \left(\hat{F}_rG_s-\hat{F}_sG_r\right)\left(\hat{H}_s\hat{K}_r-\hat{H}_r\hat{K}_s\right).
\end{align*}
This last line is \eqref{e:nNewBezout3}.  The algebra above is explained via the following steps, the second to last equality follows by combining the terms involving $\hat{H}_rH_s$ and $\hat{K}_sK_r$ by making a change of index and simple algebra. For the term involving $K_rH_s$ one makes a similar change of index argument to see that the sum from the line before is equal to a related sum, namely:
$$
\sum_{s=1}^n\sum_{r=1}^n\left(\hat{G}_{s}F_r-\hat{G}_rF_{s}\right)K_rH_s=-\sum_{s=1}^n\sum_{r=1}^n\left(\hat{G}_{s}F_r-\hat{G}_rF_{s}\right)K_sH_r;
$$
these are then combined in an obvious way to produce the claimed equality.  We further observe that the sum there is over \textit{all} $r$ and $s$, but this can be replaced by restricting to the case of $r<s$ removing the factor for $\frac{1}{2}$ since terms appear twice and the signs involved will cancel.  The term involving $\hat{H}_r\hat{K}_s$ is handled identically.

One notes that the above identities hold for any solutions to the respective equations.  In particular, if $h_j$ and $k_j$ are smooth solutions to \eqref{e:nNewBezout1} and \eqref{e:nNewBezout2}, then using \eqref{e:Fjidentity} and \eqref{e:Gjidentity} we have the following identities that will play a role later.  They are proved by either taking $\overline{\partial}$ directly or first taking $\hat{\cdot}$ and then taking $\overline{\partial}$.
 {Since the functions $F_s$ and $G_s$ are holomorphic we have:}

\begin{align}
    0 & = \sum_{r=1}^n \left(\hat{F}_rF_{s}+\hat{G}_{s}G_r\right)\overline{\partial}\hat{h}_r + \sum_{r=1}^n\left(\hat{G}_{s}F_r-\hat{G}_rF_{s}\right)\overline{\partial}k_r\quad \forall 1\leq s\leq n \label{e:dbaridentity1};\\
    0 & = \sum_{r=1}^n \left(\hat{F}_sF_{r}+\hat{G}_{r}G_s\right)\overline{\partial}h_r + \sum_{r=1}^n\left(\hat{F}_{r}G_s-\hat{F}_sG_{r}\right)\overline{\partial}\hat{k}_r\quad \forall 1\leq s\leq n \label{e:dbaridentity2};\\
    0 & =\sum_{r=1}^n \left(\hat{F}_rG_{s}-\hat{F}_{s}G_r\right)\overline{\partial}\hat{h}_r - \sum_{r=1}^n\left(\hat{F}_{s}F_r+\hat{G}_rG_{s}\right)\overline{\partial}k_r\quad \forall 1\leq s\leq n \label{e:dbaridentity3};\\
    0 & =\sum_{r=1}^n \left(\hat{G}_sF_{r}-\hat{G}_{r}F_s\right)\overline{\partial}h_r - \sum_{r=1}^n\left(\hat{F}_{r}F_s+\hat{G}_sG_{r}\right)\overline{\partial}\hat{k}_r\quad \forall 1\leq s\leq n \label{e:dbaridentity4}.
\end{align}


\subsubsection{Step 2: Constructing Smooth and Bounded Solutions}

We now show how to solve \eqref{e:nNewBezout1} and \eqref{e:nNewBezout2} in a smooth and bounded manner.  First, define:
\begin{align*}
    \varphi_j(z) & := \frac{\overline{F_j(z)}}{D(z)}, \quad \psi_j(z)  := -\frac{G_j(\overline{z})}{D(\overline{z})},\\
    D(z) & := \sum_{j=1}^n\left\vert F_j(z)\right\vert^2+\sum_{j=1}^n\left\vert G_j(z)\right\vert^2.
\end{align*}
Note that $0<\delta^2\leq D(z)\leq 1$ by the hypothesis that the $F_j$'s and $G_j$'s do not have a common zero.  It is immediate to demonstrate that
$$
\sum_{j=1}^n F_j\varphi_j-\sum_{j=1}^n G_j\hat{\psi}_j=1.
$$
This follows by direct substitution of the expressions and performing the obvious operations to arrive at an expression of the form $\frac{D(z)}{D(z)}=1$.  In particular, we have solutions with the estimates
$$
\Vert \varphi_j\Vert_{L^\infty(\mathbb{D})}\leq\frac{1}{\delta^2}\ \ {\rm and}\ \ \Vert \psi_j\Vert_{L^\infty(\mathbb{D})}\leq\frac{1}{\delta^2}\ \ {\rm for} \ \ 1\leq j\leq n.
$$

Next define

\begin{align}
h_j &:= \varphi_j+\sum_{\overset{r=1}{r\neq j}}^n\lambda_{r,j} F_r+\sum_{r=1}^n\mu_{r,j}G_r,\quad 1\leq j\leq n\label{e:hjfunctions}\\
k_j &:= \psi_j+\sum_{\overset{r=1}{r\neq j}}^n\hat{\alpha}_{r,j}\hat{G}_r+\sum_{r=1}^n\hat{\mu}_{j,r}\hat{F}_r,\quad 1\leq j\leq n\label{e:kjfunctions}.
\end{align}
where $\lambda_{k,j}=-\lambda_{j,k}$, $\alpha_{k,j}=-\alpha_{j,k}$ for all $j$ and $k$ and $\lambda_{j,k}, \alpha_{j,k},\mu_{k,j}\in L^\infty(\mathbb{D})$.  Observe that we are choosing $n^2$ functions $\mu_{k,j}$ and $\frac{n^2-n}{2}$ functions $\lambda_{k,j}$ and $\alpha_{k,j}$. Because of the conditions on $\lambda_{k,j}$, $\alpha_{k,j}$, and $\mu_{k,j}$ one checks that for any choice of these functions the functions $h_j$ and $k_j$ satisfy \eqref{e:nNewBezout1}.  Indeed, we have that
\begin{align*}
    \sum_{j=1}^n F_j h_j& -\sum_{j=1}^n G_j\hat{K}_j
    \\
    &=\sum_{j=1}^nF_j\left(\varphi_j+\sum_{\overset{r=1}{r\neq j}}^n\lambda_{r,j} F_r+\sum_{r=1}^n\mu_{r,j}G_r\right)
     -\sum_{j=1}^nG_j\left(\hat{\psi}_j+\sum_{\overset{r=1}{r\neq j}}^n\alpha_{r,j}G_k+\sum_{r=1}^n\mu_{j,r}F_r\right)\\
     & = \sum_{j=1}^nF_j\varphi_j-\sum_{j=1}^nG_j\hat{\psi}_j+\sum_{j=1}^n\sum_{\overset{r=1}{r\neq j}}\lambda_{r,j}F_jF_r+ \sum_{j=1}^n\sum_{\overset{r=1}{r\neq j}}\alpha_{r,j}G_jG_r\\
     & + \sum_{j=1}^n\sum_{r=1}^n\mu_{r,j}F_jG_r-\sum_{j=1}^n\sum_{r=1}^n\mu_{j,r}F_rG_j\\
     & = 1.
\end{align*}

With the last inequality following since $\varphi_j$ and $\psi_j$ satisfy \eqref{e:nNewBezout1}, the alternating property of $\lambda_{j,k}$ and $\alpha_{j,k}$, and a simple change of index argument (equivalently, seeing it as the inner product between a matrix $\mu$ acting on a vector $F$ paired against a vector $G$, and then the other term is simply an adjoint expression).\\

Now it suffices to choose certain specific functions $\lambda_{k,j}$, $\alpha_{k,j}$ and $\mu_{k,j}$ to satisfy \eqref{e:nNewBezout2}.  These will become obvious once we substitute the choices of $h_j$ and $k_j$ into \eqref{e:nNewBezout2}.  Doing this we see
\begin{align*}
\sum_{j=1}^n F_jk_j&+\sum_{j=1}^n G_j\hat{h}_j
\\
&
= \sum_{j=1}^nF_j\left(\psi_j+\sum_{\overset{r=1}{r\neq j}}^n\hat{\alpha}_{r,j}\hat{G}_r+\sum_{r=1}^n\hat{\mu}_{j,r}\hat{F}_r\right)+\sum_{j=1}^nG_j\left(\hat{\varphi}_j+\sum_{\overset{r=1}{r\neq j}}^n\hat{\lambda}_{r,j} \hat{F}_r+\sum_{r=1}^n\hat{\mu}_{r,j}\hat{G}_r\right)\\
& = \sum_{j=1}^n\left(F_j\psi_j+G_j\hat{\varphi}_j\right)+\sum_{j=1}^n\sum_{\overset{r=1}{r\neq j}}^n\hat{\alpha}_{r,j}\hat{G}_rF_j+\sum_{j=1}^n\sum_{\overset{r=1}{r\neq j}}^n\hat{\lambda}_{r,j} \hat{F}_rG_j\\
& + \sum_{j=1}^n\sum_{r=1}^n\hat{\mu}_{j,r}F_j\hat{F}_r+ \sum_{j=1}^n\sum_{r=1}^n\hat{\mu}_{r,j}G_j\hat{G}_r\\
& = \tau+\frac{1}{2}\sum_{j=1}^n\sum_{\overset{r=1}{r\neq j}}^n\hat{\alpha}_{r,j}\left(\hat{G}_rF_j-\hat{G}_jF_r\right)+\frac{1}{2}\sum_{j=1}^n\sum_{\overset{r=1}{r\neq j}}^n\hat{\lambda}_{r,j}\left(\hat{F}_rG_j-\hat{F}_jG_r\right)
\\
&
 {+ \sum_{j=1}^n\sum_{r=1}^n\hat{\mu}_{j,r}\left(F_j\hat{F}_r+\hat{G}_jG_r\right)}
\end{align*}
where we have used that
$$
\tau=\sum_{j=1}^n\left(F_j\psi_j+G_j\hat{\varphi}_j\right)
$$
and then made change of index arguments to change the resulting sums.  Observe that the function $\tau$ is smooth and bounded and satisfies the estimate $\Vert \tau\Vert_{L^\infty(\mathbb{D})}\lesssim \frac{\sqrt{n}}{\delta^2}$.

From this computation, it is possible to have \eqref{e:nNewBezout2} hold if and only if it is possible to choose functions $\lambda_{k,j}$, $\alpha_{k,j}$ and $\mu_{k,j}$ that are smooth and bounded for which
\begin{equation}
\label{e:nTau}
-\tau=\frac{1}{2}\sum_{j=1}^n\sum_{\overset{r=1}{r\neq j}}^n\hat{\alpha}_{r,j}\left(\hat{G}_rF_j-\hat{G}_jF_r\right)
+\frac{1}{2}\sum_{j=1}^n\sum_{\overset{r=1}{r\neq j}}^n\hat{\lambda}_{r,j}\left(\hat{F}_rG_j-\hat{F}_jG_r\right)
 {+ \sum_{j=1}^n\sum_{r=1}^n\hat{\mu}_{j,r}\left(F_j\hat{F}_r+\hat{G}_jG_r\right)}.
\end{equation}
To find these functions set:
$$
\Delta(z):=\sum_{r=1}^{n}\sum_{j=1}^n
\abs{ {F_j\hat{F}_r+G_r\hat{G}_j}}^2
+\sum_{r=1}^{n}\sum_{j=r+1}^n\abs{\hat{G}_rF_j-\hat{G}_jF_r}^2
+\sum_{r=1}^{n}\sum_{j=r+1}^n\abs{\hat{F}_rG_j-\hat{F}_jG_r}^2.
$$
One observes that \eqref{e:nNewBezout3} implies that $\Delta(z)\geq\delta^2>0$ for all $z\in\mathbb{D}$; and in particular the argument given shows that this condition is necessary for the possible solutions to Theorem \ref{t:NewCoronaS2} in the case of arbitrary $n$ and so \eqref{e:NewCC2} is a necessary reasonable condition to impose in the course of the proof.

Define
\begin{align*}
    \hat{\mu}_{k,j} & :=
    -\tau\frac{\overline{ {F_j\hat{F}_k+G_j\hat{G}_k}}}{\Delta}\quad 1\leq k\leq n, 1\leq j\leq n\\
    \hat{\alpha}_{k,j} & := -\tau\frac{\overline{\hat{G}_kF_j-\hat{G}_jF_k}}{\Delta}\quad 1\leq k\leq n, 1\leq j\leq n\\
    \hat{\lambda}_{k,j} & := -\tau\frac{\overline{\hat{F}_kG_j-\hat{F}_jG_k}}{\Delta}\quad 1\leq k\leq n, 1\leq j\leq n.
\end{align*}
By the properties of $\Delta$ and $\tau$ one observes that $\mu_{k,j},\alpha_{k,j},\lambda_{k,j}\in L^\infty(\mathbb{D})$ with norm controlled by $\frac{\sqrt{n}}{\delta^4}$.  And further that the alternating property $\alpha_{k,j}=-\alpha_{j,k}$ and $\lambda_{k,j}=-\lambda_{j,k}$ holds.  Hence, with these choices of $\mu_{k,j},\alpha_{k,j}$ and $\lambda_{k,j}$ that we have \eqref{e:nTau} holding.

This whole argument has produced smooth and bounded functions $h_j$ and $k_j$ such that \eqref{e:nNewBezout1} and \eqref{e:nNewBezout2} hold:
\begin{align*}
    \sum_{j=1}^n F_jh_j-\sum_{j=1}^nG_j\hat{k}_j & =1\\
    \sum_{j=1}^n F_jk_j+\sum_{j=1}^nG_j\hat{h}_j & =0
\end{align*}
and the following estimates:
$$
\Vert h_j\Vert_{L^\infty(\mathbb{D})}\lesssim \frac{1}{\delta^2}+\frac{n}{\delta^4}\ \ {\rm and}\ \  \Vert k_j\Vert_{L^\infty(\mathbb{D})}\lesssim \frac{1}{\delta^2}+\frac{n}{\delta^4}\ \ {\rm for}\ \  1\leq j\leq n.
$$

\subsubsection{Step 3: Perturbing the Solutions to be Holomorphic}

It remains to correct these functions $h_j$ and $k_j$ in a manner that will make holomorphic solutions.  Define
\begin{align}
    H_j & = h_j+\sum_{r=1}^{n}\sum_{\overset{s=1}{s\neq j}}^n \left(\hat{F}_rF_s+\hat{G}_sG_r\right)\beta_{r,s,j}+\sum_{r=1}^n\sum_{s=r+1}^n\left(\hat{F}_sG_r-\hat{F}_rG_s\right)\gamma_{r,s,j}\label{e:nHfunctions}\\
    & +\sum_{\overset{r=1}{r\neq j}}^n\sum_{\overset{s=r+1}{s\neq j}}^n\left(\hat{G}_rF_s-\hat{G}_sF_r\right)\eta_{r,s,j}\notag\\
    K_j & = k_j+\sum_{r=1}^{n}\sum_{\overset{s=1}{s\neq j}}^n \left(\hat{F}_rF_s+\hat{G}_sG_r\right)\hat{\gamma}_{j,s,r}+\sum_{r=1}^n\sum_{s=r+1}^n\left(\hat{F}_rG_s-\hat{F}_sG_r\right)\hat{\beta}_{j,s,r}\label{e:nKfunctions}\\
    & +\sum_{\overset{r=1}{r\neq j}}^n\sum_{\overset{s=r+1}{s\neq j}}^n
    \left(\hat{G}_rF_s-\hat{G}_sF_r\right)\hat{\tilde{\eta}}_{r,s,j}\notag
\end{align}
where the functions $\beta_{r,s,j}$ satisfy $\beta_{r,s,j}=-\beta_{r,j,s}$ for $1\leq r\leq n$, $1\leq s\leq n$, and $1\leq j\leq n$; the functions $\gamma_{r,s,j}$ satisfy the alternating condition $\gamma_{r,s,j}=-\gamma_{s,r,j}$ for $1\leq r,s\leq n$.  Note that the alternating condition in the $\beta_{r,s,j}$ is in the second and third indices, while the alternating condition for the $\gamma_{r,s,j}$ is in the first two indices!  We impose the condition that $\eta_{r,s,j}=-\eta_{j,s,r}=\eta_{j,r,s}$ for $1\leq r\leq n$, $1\leq s\leq n$ and $1\leq j\leq n$; an analogous condition is imposed on the functions $\tilde{\eta}_{r,s,j}$.

With these functions defined it is the case that they solve \eqref{e:nNewBezout1}.  Substituting in and performing obvious computations and collecting terms we have:
\begin{align*}
    \sum_{j=1}^n F_jH_j& -\sum_{j=1}^nG_j\hat{K}_j
     \\
     &
     = \sum_{j=1}^n F_jh_j+\sum_{j=1}^n F_j\sum_{r=1}^{n}\sum_{\overset{s=1}{s\neq j}}^n \left(\hat{F}_rF_s+\hat{G}_sG_r\right)\beta_{r,s,j}\\
    & +\sum_{j=1}^n F_j\sum_{r=1}^n\sum_{s=r+1}^n\left(\hat{F}_sG_r-\hat{F}_rG_s\right)\gamma_{r,s,j}+\sum_{j=1}^nF_j\sum_{\overset{r=1}{r\neq j}}^n\sum_{\overset{s=r+1}{s\neq j}}^n\left(\hat{G}_rF_s-\hat{G}_sF_r\right)\eta_{r,s,j}\\
    & -\sum_{j=1}^nG_j\hat{k}_j-\sum_{j=1}^nG_j\sum_{r=1}^{n}\sum_{\overset{s=1}{s\neq j}}^n \left(\hat{F}_sF_r+\hat{G}_rG_s\right)\gamma_{j,s,r}
    \end{align*}
\begin{align*}
    & -\sum_{j=1}^nG_j\sum_{r=1}^n\sum_{s=r+1}^n\left(\hat{G}_sF_r-\hat{G}_rF_s\right)\beta_{j,s,r}-\sum_{j=1}^nG_j\sum_{\overset{r=1}{r\neq j}}^n\sum_{\overset{s=r+1}{s\neq j}}^n\left(\hat{F}_sG_r-\hat{F}_rG_s\right)\tilde{\eta}_{r,s,j}\\
    & = \sum_{j=1}^n \left(F_jh_j-G_j\hat{k}_j\right)+I+II+III+IV\\
    & = 1+I+II+III+IV=1
\end{align*}
with the term $I$, $II$ and $III$ and $IV$ being defined as
\begin{align}
    I & :=\sum_{j=1}^n F_j\sum_{r=1}^{n}\sum_{\overset{s=1}{s\neq j}}^n \left(\hat{F}_rF_s+\hat{G}_sG_r\right)\beta_{r,s,j}-\sum_{j=1}^nG_j\sum_{r=1}^n\sum_{s=r+1}^n\left(\hat{G}_sF_r-\hat{G}_rF_s\right)\beta_{j,s,r};\label{e:identities1}\\
    II & := \sum_{j=1}^n F_j\sum_{r=1}^n\sum_{s=r+1}^n\left(\hat{F}_sG_r-\hat{F}_rG_s\right)\gamma_{r,s,j}-\sum_{j=1}^nG_j\sum_{r=1}^{n}\sum_{\overset{s=1}{s\neq j}}^n \left(\hat{F}_sF_r+\hat{G}_rG_s\right)\gamma_{j,s,r}\label{e:identities2};\\
    III & := \sum_{j=1}^nF_j\sum_{\overset{r=1}{r\neq j}}^n\sum_{\overset{s=r+1}{s\neq j}}^n\left(\hat{G}_rF_s-\hat{G}_sF_r\right)\eta_{r,s,j}\label{e:identities3};\\
    IV & :=-\sum_{j=1}^nG_j\sum_{\overset{r=1}{r\neq j}}^n\sum_{\overset{s=r+1}{s\neq j}}^n\left(\hat{F}_sG_r-\hat{F}_rG_s\right)\tilde{\eta}_{r,s,j}\label{e:identities4};
\end{align}
and in the last line we used that $h_j$ and $k_j$ solve \eqref{e:nNewBezout1} and that $I=II=III=IV=0$, that we now turn to demonstrating.

The idea to show that these terms are zero is to write sums like $\displaystyle\sum_{s\neq j}$ as $\displaystyle\sum_{s<j}+\sum_{j<s}$ and then use a reindexing of a double sum to be over the same sets (essentially summing over columns and then rows versus rows and then columns), use a change of variables/index to arrive at similar expressions, or, use a Fubini argument to change the order of summation and the resulting summation sets.  The alternating conditions on the function $\beta_{r,s,j}$ will result in a summation over a common term times a $\beta_{r,s,j}$, which will collapse to zero.  Analogous ideas are applied to the terms involving $\gamma_{r,s,j}$ and similar ideas are used in the expressions involving $\eta_{r,s,j}$ and $\tilde{\eta}_{r,s,j}$.

We do this now for term $I$, with each equality following from ones of the steps of the type outlined above:
\begin{align*}
    I & =\sum_{j=1}^n F_j\sum_{r=1}^{n}\sum_{\overset{s=1}{s\neq j}}^n \left(\hat{F}_rF_s+\hat{G}_sG_r\right)\beta_{r,s,j}-\sum_{j=1}^nG_j\sum_{r=1}^n\sum_{s=r+1}^n\left(\hat{G}_sF_r-\hat{G}_rF_s\right)\beta_{j,s,r}\\
    & = \sum_{r=1}^n\left(\sum_{j=1}^n\sum_{s=j+1}^n+\sum_{j=1}^{n}\sum_{s=1}^{j-1}\right) \left(F_j \left(\hat{F}_rF_s+\hat{G}_sG_r\right)\beta_{r,s,j}\right)\\
    & \ \ \ \ \ \ \ -\sum_{j=1}^n\sum_{r=1}^n\sum_{s=r+1}^nG_j\left(\hat{G}_sF_r-\hat{G}_rF_s\right)\beta_{j,s,r}\\
    & = \sum_{r=1}^n\sum_{j=1}^n\sum_{s=j+1}^n\left( F_j\left(\hat{F}_rF_s+\hat{G}_sG_r\right)\beta_{r,s,j}+F_s\left(\hat{F}_rF_j+\hat{G}_jG_r\right)\beta_{r,j,s}\right)
    \\
    & \ \ \ \ \ \ \ -\sum_{r=1}^n\sum_{j=1}^n\sum_{s=j+1}^n G_r\left(\hat{G}_sF_j-\hat{G}_jF_s\right)\beta_{r,s,j}
    \\
     & =
     \sum_{r=1}^n\sum_{j=1}^n\sum_{s=j+1}^n \left(-F_j\left(\hat{F}_rF_s+\hat{G}_sG_r\right)+F_s\left(\hat{F}_rF_j+\hat{G}_jG_r\right)
    +G_r\left(\hat{G}_sF_j-\hat{G}_jF_s\right)\right)\beta_{r,j,s}
    \\
    & = \sum_{r=1}^n\sum_{j=1}^n\sum_{s=j+1}^n \left(-F_j\hat{F}_rF_s-F_j\hat{G}_sG_r+F_s\hat{F}_rF_j+F_s\hat{G}_jG_r
    +G_r\hat{G}_sF_j-G_r\hat{G}_jF_s\right)\beta_{r,j,s}=0,
\end{align*}
with the last line following by obvious algebra.  Term $II$ is handled analogously by observing that is effectively the same as $I$ but the roles of $F_j$ and $G_j$ are interchanged. We handle term $III$ next; the proof of this term will of course handle $IV$ by a symmetric argument.  Observe that:
\begin{align*}
     III & = \sum_{j=1}^nF_j\sum_{\overset{r=1}{r\neq j}}^n\sum_{\overset{s=r+1}{s\neq j}}^n\left(\hat{G}_rF_s-\hat{G}_sF_r\right)\eta_{r,s,j}\\
     &= \sum_{j=1}^n\sum_{r=j+1}^n\sum_{s=r+1}^n F_j\left(\hat{G}_rF_s-\hat{G}_sF_r\right)\eta_{r,s,j}+\sum_{j=1}^n\sum_{r=1}^{j-1}\sum_{s=r+1}^{j-1}F_j\left(\hat{G}_rF_s
     -\hat{G}_sF_r\right)\eta_{r,s,j}
     \\
     & +\sum_{j=1}^n\sum_{r=1}^{j-1}\sum_{s=j+1}^{n}F_j\left(\hat{G}_rF_s-\hat{G}_sF_r\right)\eta_{r,s,j}\\
    & = \sum_{j=1}^n\sum_{r=j+1}^n\sum_{s=r+1}^n F_j\left(\hat{G}_rF_s-\hat{G}_sF_r\right)\eta_{r,s,j}+\sum_{r=1}^n\sum_{s=r+1}^n\sum_{j=s+1}^nF_j\left(\hat{G}_rF_s-\hat{G}_sF_r\right)\eta_{r,s,j}\\
     & +\sum_{r=1}^n\sum_{j=r+1}^n\sum_{s=r+1}^nF_j\left(\hat{G}_rF_s-\hat{G}_sF_r\right)\eta_{r,s,j}
      \end{align*}
     
     \begin{align*}
     \quad & = \sum_{j=1}^n\sum_{r=j+1}^n\sum_{s=r+1}^n F_j\left(\hat{G}_rF_s-\hat{G}_sF_r\right)\eta_{r,s,j}+\sum_{j=1}^n\sum_{r=j+1}^n\sum_{s=r+1}^nF_s\left(\hat{G}_jF_r-\hat{G}_rF_j\right)\eta_{r,s,j}\\
     & -\sum_{j=1}^n\sum_{r=j+1}^n\sum_{s=r+1}^nF_r\left(\hat{G}_jF_s-\hat{G}_sF_j\right)\eta_{r,s,j}\\
     & = \sum_{j=1}^n\sum_{r=j+1}^n\sum_{s=r+1}^n \left(F_j\left(\hat{G}_rF_s-\hat{G}_sF_r\right) + F_s\left(\hat{G}_jF_r-\hat{G}_rF_j\right)-F_r\left(\hat{G}_jF_s-\hat{G}_sF_j\right) \right)\eta_{r,s,j}\\
     & = \sum_{j=1}^n\sum_{r=j+1}^n\sum_{s=r+1}^n \left(F_j\hat{G}_rF_s-F_j\hat{G}_sF_r + F_s\hat{G}_jF_r-F_s\hat{G}_rF_j-F_r\hat{G}_jF_s+F_r\hat{G}_sF_j\right)\eta_{r,s,j}
     \\
     &=0.
\end{align*}

Next, it will be the case that the functions $H_j$ and $K_j$ in \eqref{e:nHfunctions} and \eqref{e:nKfunctions} satisfy \eqref{e:nNewBezout2}.  Making the substitutions and obvious computations we have:
\begin{align*}
    \sum_{j=1}^n F_jK_j+\sum_{j=1}^n G_j\hat{H}_j & = \sum_{j=1}^n F_jk_j+\sum_{j=1}^n F_j \sum_{r=1}^{n}\sum_{\overset{s=1}{s\neq j}}^n \left(\hat{F}_rF_s+\hat{G}_sG_r\right)\hat{\gamma}_{j,s,r}\\
    & +\sum_{j=1}^nF_j\sum_{r=1}^n\sum_{s=r+1}^n
    \left(\hat{F}_rG_s-\hat{F}_sG_r\right)\hat{\beta}_{j,s,r}+\sum_{j=1}^n F_j
    \sum_{\overset{s=r+1}{s\neq j}}^n\left(\hat{G}_rF_s-\hat{G}_sF_r\right)\hat{\tilde{\eta}}_{r,s,j}\\
    & + \sum_{j=1}^nG_j\hat{h}_j + \sum_{j=1}^n G_j \sum_{r=1}^{n}\sum_{\overset{s=1}{s\neq j}}^n \left(\hat{F}_sF_r+\hat{G}_rG_s\right)\hat{\beta}_{r,s,j}\\
    & +\sum_{j=1}^nG_j\sum_{r=1}^n\sum_{s=r+1}^n\left(\hat{G}_rF_s
    -\hat{G}_sF_r\right)\hat{\gamma}_{r,s,j}+\sum_{j=1}^nG_j\sum_{\overset{s=r+1}{s\neq j}}^n\left(\hat{F}_sG_r-\hat{F}_rG_s\right)\hat{\eta}_{r,s,j}\\
    & = \sum_{j=1}^n \left(F_jk_j+G_j\hat{h}_j\right)+I+II+III+IV\\
    & = 0+I+II+III+IV=0
\end{align*}
where similar to above, we define $I$, $II$, $III$ and $IV$ as
\begin{align*}
    I & := \sum_{j=1}^n F_j \sum_{r=1}^{n}\sum_{\overset{s=1}{s\neq j}}^n \left(\hat{F}_rF_s+\hat{G}_sG_r\right)\hat{\gamma}_{j,s,r}+\sum_{j=1}^nG_j\sum_{r=1}^n\sum_{s=r+1}^n\left(\hat{G}_rF_s-\hat{G}_sF_r\right)\hat{\gamma}_{r,s,j};\\
    II & := \sum_{j=1}^nF_j\sum_{r=1}^n\sum_{s=r+1}^n\left(\hat{F}_rG_s-\hat{F}_sG_r\right)\hat{\beta}_{j,s,r}+ \sum_{j=1}^n G_j \sum_{r=1}^{n}\sum_{\overset{s=1}{s\neq j}}^n \left(\hat{F}_sF_r+\hat{G}_rG_s\right)\hat{\beta}_{r,s,j};\\
    III & := \sum_{j=1}^n F_j \sum_{\overset{r=1}{r\neq j}}^n\sum_{\overset{s=r+1}{s\neq j}}^n\left(\hat{G}_rF_s-\hat{G}_sF_r\right)\hat{\tilde{\eta}}_{r,s,j};\\
    IV & := \sum_{j=1}^n G_j\sum_{\overset{r=1}{r\neq j}}^n\sum_{\overset{s=r+1}{s\neq j}}^n\left(\hat{F}_sG_r-\hat{F}_rG_s\right)\hat{\eta}_{r,s,j}.
\end{align*}
But, observe that these terms are basically the same as terms $I$, $II$, $III$ and $IV$ in \eqref{e:identities1}-\eqref{e:identities4} and are then zero as claimed via similar computations that we omit.

Now, with the most general solution to \eqref{e:nNewBezout1} and \eqref{e:nNewBezout2} we proceed to select the $\beta_{r,s,j}$, $\gamma_{r,s,j}$, $\eta_{r,s,j}$ and $\tilde{\eta}_{r,s,j}$ to satisfy certain $\overline{\partial}$-equations.  In particular, we demand that 
\begin{align}
    \overline{\partial}\beta_{r,s,j} & := \overline{\partial}h_s\left(\hat{h}_rh_j+\hat{k}_rk_j\right)-\overline{\partial}h_j\left(\hat{h}_rh_s+\hat{k}_rk_s\right)-\overline{\partial}\hat{k}_r\left(h_sk_j-h_jk_s\right)\label{e:dbetarsj}\\
    \overline{\partial}\gamma_{r,s,j} & := \overline{\partial}\hat{k}_s\left(\hat{h}_rh_j+\hat{k}_rk_j\right)-\overline{\partial}\hat{k}_r\left(\hat{h}_sh_j+\hat{k}_sk_j\right)-\overline{\partial}h_j\left(\hat{h}_r\hat{k}_s-\hat{h}_s\hat{k}_r\right)\label{e:dgammarsj}\\
    \overline{\partial}\eta_{r,s,j} & :=- \overline{\partial}h_s\left(h_jk_r-h_rk_j\right)-\overline{\partial}h_r\left(h_sk_j-h_jk_s\right)-\overline{\partial}h_j\left(h_rk_s-h_sk_r\right)\label{e:detarsj}\\
    \overline{\partial}\tilde{\eta}_{r,s,j} & :=- \overline{\partial}\hat{k}_s\left(\hat{h}_r\hat{k}_j-\hat{h}_j\hat{k}_r\right)-\overline{\partial}\hat{k}_r\left(\hat{h}_j\hat{k}_s-\hat{h}_s\hat{k}_j\right)-\overline{\partial}\hat{k}_j\left(\hat{h}_s\hat{k}_r-\hat{h}_r\hat{k}_s\right)\label{e:detatildersj}.
\end{align}
\begin{rem}
    The choice of $\beta_{r,s, j}$, $\gamma_{r,s,j}$, $\eta_{r,s,J}$ and $\tilde{\eta}_{r,s,j}$ arise from taking the negative of determinants of $3\times 3$ minors of the following matrix:
    $$
    \left(\begin{array}{cccccccccc}
      h_1 & h_2 & \cdots & h_{n-1} & h_n & -\hat{k}_1 & -\hat{k}_2 & \cdots & -\hat{k}_{n-1} & -\hat{k}_n\\
      k_1 & k_2 & \cdots & k_{n-1} & k_n & \hat{h}_1 & \hat{h}_2 & \cdots & \hat{h}_{n-1} & \hat{h}_n\\
      \overline{\partial}h_1 & \overline{\partial}h_2 & \cdots & \overline{\partial}h_{n-1} & \overline{\partial}h_n & -\overline{\partial}\hat{k}_1 & - \overline{\partial}\hat{k}_2 & \cdots & -\overline{\partial}\hat{k}_{n-1} & -\overline{\partial}\hat{k}_n
    \end{array}\right).
    $$
    In particular, the $\beta_{r,s,j}$ arises from taking two columns of the $h_j$ terms and one column of the $\hat{k}_j$; $\gamma_{r,s,j}$ arises from one column of the $h_j$ and two columns of the $\hat{k}_j$; $\eta_{r,s,j}$ from taking three columns of the $h_j$ terms; and $\tilde{\eta}_{r,s,j}$ from taking three columns of the $\hat{k}_j$ terms.  By elementary properties of determinants it holds that $\beta_{r,s,j}=-\beta_{r,j,s}$, $\gamma_{r,s,j}=-\gamma_{s,r,j}$, $\eta_{r,s,j}=-\eta_{r,j,s}=\eta_{j,r,s}$ and $\tilde{\eta}_{r,s,j}=-\tilde{\eta}_{r,j,s}=\tilde{\eta}_{j,r,s}$ as needed.
\end{rem}

With this choice of $\beta_{r,s,j}$, $\gamma_{r,s,j}$, $\eta_{r,s,j}$ and $\tilde{\eta}_{r,s,j}$ the functions $H_j$ and $K_j$ are holomorphic.  Observe that from \eqref{e:nHfunctions} and \eqref{e:nKfunctions} that:
\begin{align*}
    \overline{\partial}H_j & = \overline{\partial}h_j+\sum_{r=1}^{n}\sum_{\overset{s=1}{s\neq j}}^n \left(\hat{F}_rF_s+\hat{G}_sG_r\right)\overline{\partial}\beta_{r,s,j}+\sum_{r=1}^n\sum_{s=r+1}^n\left(\hat{F}_sG_r-\hat{F}_rG_s\right)\overline{\partial}\gamma_{r,s,j}\\
    & +\sum_{\overset{r=1}{r\neq j}}^n\sum_{\overset{s=r+1}{s\neq j}}^n\left(\hat{G}_rF_s-\hat{G}_sF_r\right)\overline{\partial}\eta_{r,s,j};\\
    \overline{\partial}K_j & = \overline{\partial}k_j+\sum_{r=1}^{n}\sum_{\overset{s=1}{s\neq j}}^n \left(\hat{F}_rF_s+\hat{G}_sG_r\right)\overline{\partial}\hat{\gamma}_{j,s,r}+\sum_{r=1}^n\sum_{s=r+1}^n\left(\hat{F}_rG_s-\hat{F}_sG_r\right)\overline{\partial}\hat{\beta}_{j,s,r}\\
    & +\sum_{\overset{r=1}{r\neq j}}^n\sum_{\overset{s=r+1}{s\neq j}}^n
    \left(\hat{G}_rF_s-\hat{G}_sF_r\right)\overline{\partial}\hat{\tilde{\eta}}_{r,s,j}.
\end{align*}

We claim that for $1\leq j\leq n$,
\begin{align}
    -\overline{\partial}h_j & = \sum_{r=1}^{n}\sum_{\overset{s=1}{s\neq j}}^{n} \left(\hat{F}_rF_s+\hat{G}_sG_r\right)\overline{\partial}\beta_{r,s,j} +\sum_{r=1}^{n}\sum_{s=r+1}^{n}\left(\hat{F}_sG_r-\hat{F}_rG_s\right)\overline{\partial}\gamma_{r,s,j}\label{e:dhjmagic}\\
     & +\sum_{\overset{r=1}{r\neq j}}^{n}\sum_{\overset{s=r+1}{s\neq j}}^{n}\left(\hat{G}_rF_s-\hat{G}_sF_r\right)\overline{\partial}\eta_{r,s,j};\notag\\
    -\overline{\partial}k_j & = \sum_{r=1}^{n}\sum_{\overset{s=1}{s\neq j}}^{n} \left(\hat{F}_rF_s+\hat{G}_sG_r\right)\overline{\partial}\hat{\gamma}_{j,s,r} +\sum_{r=1}^{n}\sum_{s=r+1}^{n}\left(\hat{F}_rG_s-\hat{F}_sG_r\right)\overline{\partial}\hat{\beta}_{j,s,r}\label{e:dkjmagic}\\
     & +\sum_{\overset{r=1}{r\neq j}}^{n}\sum_{\overset{s=r+1}{s\neq j}}^{n}\left(\hat{G}_rF_s-\hat{G}_sF_r\right)\overline{\partial}\hat{\tilde{\eta}}_{r,s,j}\notag
\end{align}
which of course then implies that $\overline{\partial}H_j =0$ and $\overline{\partial}K_j=0$ are holomorphic solutions to \eqref{e:nNewBezout1} and \eqref{e:nNewBezout2} and the functions in \eqref{e:nHfunctions} and \eqref{e:nKfunctions} are holomorphic.  We provide the argument for \eqref{e:dhjmagic} as the proof of \eqref{e:dkjmagic} is done similarly.  Except in this case, we use the observation that $\overline{\partial}\hat{\varphi}=\widehat{\overline{\partial}\varphi}$.  This follows from the chain rule for derivatives involving $\overline{\partial}$ and $\partial$ and basic properties of $\overline{\partial}$ and $\partial$.  Indeed, one sees that
$$
\overline{\partial}\hat{\varphi} (z)= \partial \overline{\varphi}(\overline{z})\overline{\partial}\overline{z}+\overline{\partial}\varphi(\overline{z})\overline{\partial}z=\partial\overline{\varphi}(\overline{z})=\overline{\overline{\partial}\varphi(\overline{z})}=\widehat{\overline{\partial}\varphi}(z).
$$

Turning to the proof of \eqref{e:dhjmagic}, first, observe that by \eqref{e:nNewBezout3} we have the following holds:
\begin{align*}
    -\overline{\partial}h_j & = -\sum_{s=1}^{n}\sum_{r=1}^n \left(\hat{F}_rF_s+\hat{G}_sG_r\right)\left(\hat{h}_rh_s+\hat{k}_rk_s\right)\overline{\partial}h_j
    \\
    &
    -\sum_{r=1}^{n}\sum_{s=r+1}^n \left(\hat{G}_rF_s-\hat{G}_sF_r\right)\left(k_sh_r-k_rh_s\right)\overline{\partial}h_j\\
    & -\sum_{r=1}^{n}\sum_{s=r+1}^n \left(\hat{F}_sG_r-\hat{F}_rG_s\right)\left(\hat{h}_r\hat{k}_s-\hat{h}_s\hat{k}_r\right)\overline{\partial}h_j\\
    & := T_1+T_2+T_3
\end{align*}
where $T_1$, $T_2$ and $T_3$ are defined by the expressions above.

Expand $T_1$, $T_2$ and $T_3$ via the formulas \eqref{e:dbetarsj}, \eqref{e:dgammarsj} and \eqref{e:detarsj}, and then utilize identities \eqref{e:dbaridentity2} and \eqref{e:dbaridentity4}.  Doing this, yields the following identities:
\begin{align}
    T_1 & = \sum_{r=1}^{n}\sum_{\overset{s=1}{s\neq j}}^{n} \left(\hat{F}_rF_s+\hat{G}_sG_r\right)\overline{\partial}\beta_{r,s,j}+\sum_{r=1}^{n}\sum_{s=1}^{n} \left(\hat{F}_sG_r-\hat{F}_rG_s\right)\left(\hat{h}_rh_j+\hat{k}_rk_j\right)\overline{\partial}\hat{k}_s\label{e:T1}\\
    & +\sum_{r=1}^{n}\sum_{s=1}^{n}\left(\hat{G}_sF_r-\hat{G}_rF_s\right)\left(h_sk_j-h_jk_s\right)\overline{\partial}h_r = \sum_{r=1}^{n}\sum_{\overset{s=1}{s\neq j}}^{n}\left(\hat{F}_rF_s+\hat{G}_sG_r\right)\overline{\partial}\beta_{r,s,j}+T_{1,1}+T_{1,2}\notag
\end{align}
where
\begin{align*}
T_{1,1} & :=\sum_{r=1}^{n}\sum_{s=1}^{n} \left(\hat{F}_sG_r-\hat{F}_rG_s\right)\left(\hat{h}_rh_j+\hat{k}_rk_j\right)\overline{\partial}\hat{k}_s\\
T_{1,2} & := \sum_{r=1}^{n}\sum_{s=1}^{n}\left(\hat{G}_sF_r-\hat{G}_rF_s\right)\left(h_sk_j-h_jk_s\right)\overline{\partial}h_r.
\end{align*}
And
\begin{align}
    T_2 & = \sum_{\overset{r=1}{r\neq j}}^{n}\sum_{\overset{s=r+1}{s\neq j}}^{n}\left(\hat{G}_rF_s-\hat{G}_sF_r\right)\overline{\partial}\eta_{r,s,j}+\sum_{r=1}^n\sum_{s=r+1}^{n}\left(\hat{G}_rF_s-\hat{G}_sF_r\right)\left(h_jk_r-h_rk_j\right)\overline{\partial}h_s \label{e:T2}\\
    & +\sum_{r=1}^n\sum_{s=r+1}^{n}\left(\hat{G}_rF_s-\hat{G}_sF_r\right)\left(h_sk_j-h_jk_s\right)\overline{\partial}h_r\notag
     \\
     &
     = \sum_{\overset{r=1}{r\neq j}}^{n}\sum_{\overset{s=r+1}{s\neq j}}^{n}\left(\hat{G}_rF_s-\hat{G}_sF_r\right)\overline{\partial}\eta_{r,s,j} + T_{2,1}+T_{2,2}\notag
\end{align}
where
\begin{align*}
    T_{2,1} := \sum_{r=1}^n\sum_{s=r+1}^{n}\left(\hat{G}_rF_s-\hat{G}_sF_r\right)\left(h_jk_r-h_rk_j\right)\overline{\partial}h_s\\
    T_{2,2} := \sum_{r=1}^n\sum_{s=r+1}^{n}\left(\hat{G}_rF_s-\hat{G}_sF_r\right)\left(h_sk_j-h_jk_s\right)\overline{\partial}h_r.
\end{align*}
Finally,
\begin{align}
    T_3 & = \sum_{r=1}^{n}\sum_{s=r+1}^{n}\left(\hat{F}_sG_r-\hat{F}_rG_s\right)\overline{\partial}\gamma_{r,s,j}-\sum_{r=1}^n\sum_{s=r+1}^{n}\left(\hat{F}_sG_r-\hat{F}_rG_s\right)\left(\hat{h}_rh_j+\hat{k}_rk_j\right)\overline{\partial}\hat{k}_s\label{e:T3}\\
    & +\sum_{r=1}^n\sum_{s=r+1}^{n}\left(\hat{F}_sG_r-\hat{F}_rG_s\right)\left(\hat{h}_sh_j+\hat{k}_sk_j\right)\overline{\partial}\hat{k}_r
     \notag\\
     &
     = \sum_{r=1}^{n}\sum_{s=r+1}^{n}\left(\hat{F}_sG_r-\hat{F}_rG_s\right)\overline{\partial}\gamma_{r,s,j}+T_{3,1}+T_{3,2}\notag
\end{align}
with
\begin{align*}
    T_{3,1} & := -\sum_{r=1}^n\sum_{s=r+1}^{n}\left(\hat{F}_sG_r-\hat{F}_rG_s\right)\left(\hat{h}_rh_j+\hat{k}_rk_j\right)\overline{\partial}\hat{k}_s\\
    T_{3,2} & := \sum_{r=1}^n\sum_{s=r+1}^{n}\left(\hat{F}_sG_r-\hat{F}_rG_s\right)\left(\hat{h}_sh_j+\hat{k}_sk_j\right)\overline{\partial}\hat{k}_r.\\
\end{align*}
It is the case that
\begin{align*}
    T_{1,1}+T_{3,1}+T_{3,2} =0\quad\textnormal{and} \quad
    T_{1,2}+T_{2,1}+T_{2,2} =0.
\end{align*}
We show that the first term satisfies the claim; the other follows by analogous computations.  Indeed, we have that:
\begin{align*}
&T_{1,1}+T_{3,1}+T_{3,2}
\\
&
:= \sum_{r=1}^{n}\sum_{s=1}^{n} \left(\hat{F}_sG_r-\hat{F}_rG_s\right)\left(\hat{h}_rh_j+\hat{k}_rk_j\right)\overline{\partial}\hat{k}_s - \sum_{r=1}^n\sum_{s=r+1}^{n}\left(\hat{F}_sG_r-\hat{F}_rG_s\right)\left(\hat{h}_rh_j+\hat{k}_rk_j\right)\overline{\partial}\hat{k}_s\\
& +\sum_{r=1}^n\sum_{s=r+1}^{n}\left(\hat{F}_sG_r-\hat{F}_rG_s\right)\left(\hat{h}_sh_j+\hat{k}_sk_j\right)\overline{\partial}\hat{k}_r\\
& = \sum_{r=1}^{n}\sum_{s=r+1}^{n} \left(\hat{F}_sG_r-\hat{F}_rG_s\right)\left(\hat{h}_rh_j+\hat{k}_rk_j\right)\overline{\partial}\hat{k}_s+ \sum_{r=1}^{n}\sum_{s=1}^{r-1} \left(\hat{F}_sG_r-\hat{F}_rG_s\right)\left(\hat{h}_rh_j+\hat{k}_rk_j\right)\overline{\partial}\hat{k}_s\\
& - \sum_{r=1}^n\sum_{s=r+1}^{n}\left(\hat{F}_sG_r-\hat{F}_rG_s\right)\left(\hat{h}_rh_j+\hat{k}_rk_j\right)\overline{\partial}\hat{k}_s +\sum_{r=1}^n\sum_{s=r+1}^{n}\left(\hat{F}_sG_r-\hat{F}_rG_s\right)\left(\hat{h}_sh_j+\hat{k}_sk_j\right)\overline{\partial}\hat{k}_r\\
& = -\sum_{r=1}^{n}\sum_{s=1}^{r-1} \left(\hat{F}_sG_r-\hat{F}_rG_s\right)\left(\hat{h}_rh_j+\hat{k}_rk_j\right)\overline{\partial}\hat{k}_s -\sum_{r=1}^n\sum_{s=r+1}^{n}\left(\hat{F}_sG_r-\hat{F}_rG_s\right)\left(\hat{h}_sh_j+\hat{k}_sk_j\right)\overline{\partial}\hat{k}_r\\
& = \sum_{s=1}^{n}\sum_{r=s+1}^{n} \left(\hat{F}_sG_r-\hat{F}_rG_s\right)\left(\hat{h}_rh_j+\hat{k}_rk_j\right)\overline{\partial}\hat{k}_s +\sum_{r=1}^n\sum_{s=r+1}^{n}\left(\hat{F}_sG_r-\hat{F}_rG_s\right)\left(\hat{h}_sh_j+\hat{k}_sk_j\right)\overline{\partial}\hat{k}_r\\
& = -\sum_{r=1}^{n}\sum_{s=r+1}^{n} \left(\hat{F}_sG_r-\hat{F}_rG_s\right)\left(\hat{h}_sh_j+\hat{k}_sk_j\right)\overline{\partial}\hat{k}_r +\sum_{r=1}^n\sum_{s=r+1}^{n}\left(\hat{F}_sG_r-\hat{F}_rG_s\right)\left(\hat{h}_sh_j+\hat{k}_sk_j\right)\overline{\partial}\hat{k}_r\\
& =0;
\end{align*}
as explanation of the algebra here, the first equality is from the definition of the quantities $T_{1,1}$, $T_{3,1}$ and $T_{3,2}$, the second equality splits $T_{1,1}$ into the upper and lower portions, the third equality follows from canceling a common term, the fourth equality follows from interchanging the sum in $r$ and $s$, and then the second to last line will follow by exchanging indices, and the last is then immediate.

Using this we see that:
\begin{align*}
    T_1+T_2+T_3 & =\sum_{r=1}^{n}\sum_{\overset{s=r+1}{s\neq j}}^{n} \left(\hat{F}_rF_s+\hat{G}_sG_r\right)\overline{\partial}\beta_{r,s,j}+ \sum_{r=1}^{n}\sum_{s=r+1}^{n}\left(\hat{F}_sG_r-\hat{F}_rG_s\right)\overline{\partial}\gamma_{r,s,j}\\
    & +\sum_{\overset{r=1}{r\neq j}}^{n}\sum_{\overset{s=r+1}{s\neq j}}^{n}\left(\hat{G}_rF_s-\hat{G}_sF_r\right)\overline{\partial}\eta_{r,s,j}
\end{align*}
and hence we have
\begin{align*}
    -\overline{\partial}h_j &= \sum_{r=1}^{n}\sum_{\overset{s=1}{s\neq j}}^{n} \left(\hat{F}_rF_s+\hat{G}_sG_r\right)\overline{\partial}\beta_{r,s,j} +\sum_{r=1}^{n}\sum_{s=r+1}^{n}\left(\hat{F}_sG_r-\hat{F}_rG_s\right)\overline{\partial}\gamma_{r,s,j}
     \\
     &
     +\sum_{\overset{r=1}{r\neq j}}^{n}\sum_{\overset{s=r+1}{s\neq j}}^{n}\left(\hat{G}_rF_s
    -\hat{G}_sF_r\right)\overline{\partial}\eta_{r,s,j}
\end{align*}
which is nothing other than \eqref{e:dhjmagic}.

We now prove \eqref{e:T1}:
\begin{align*}
    T_1 & := -\sum_{s=1}^{n}\sum_{r=1}^n \left(\hat{F}_rF_s+\hat{G}_sG_r\right)\left(\hat{h}_rh_s+\hat{k}_rk_s\right)\overline{\partial}h_j\\
    & =\sum_{s=1}^{n}\sum_{r=1}^n \left(\hat{F}_rF_s+\hat{G}_sG_r\right)\left(\overline{\partial}\beta_{r,s,j}-\left(\hat{h}_rh_j+\hat{k}_rk_j\right)\overline{\partial}h_s+\left(h_sk_j-h_jk_s\right)\overline{\partial}\hat{k}_r\right)\\
    & = \sum_{s=1}^{n}\sum_{r=1}^n \left(\hat{F}_rF_s+\hat{G}_sG_r\right)\overline{\partial}\beta_{r,s,j}-\sum_{s=1}^{n}\sum_{r=1}^n \left(\hat{F}_rF_s+\hat{G}_sG_r\right)\left(\hat{h}_rh_j+\hat{k}_rk_j\right)\overline{\partial}h_s\\
    & + \sum_{s=1}^{n}\sum_{r=1}^n \left(\hat{F}_rF_s+\hat{G}_sG_r\right)\left(h_sk_j-h_jk_s\right)\overline{\partial}\hat{k}_r\\
    & = \sum_{r=1}^n \sum_{\overset{s=1}{s\neq j}}^{n} \left(\hat{F}_rF_s+\hat{G}_sG_r\right)\overline{\partial}\beta_{r,s,j}+\sum_{s=1}^{n}\sum_{r=1}^n \left(\hat{F}_sG_r-\hat{F}_rG_s\right)\left(\hat{h}_rh_j+\hat{k}_rk_j\right)\overline{\partial}\hat{h}_s\\
    & + \sum_{s=1}^{n}\sum_{r=1}^n \left(\hat{G}_sF_r-\hat{G}_rF_s\right)\left(h_sk_j-h_jk_s\right)\overline{\partial}h_r.
\end{align*}
The steps are justified in the following manner, the first equality is definition, the second uses \eqref{e:dbetarsj}, the third expands the algebra, the fourth uses the alternating property of the $\beta_{r,s,j}$ and utilize identities \eqref{e:dbaridentity2} and \eqref{e:dbaridentity4}.  

The proofs of \eqref{e:T2} and \eqref{e:T3} are done via similar arguments.

\begin{rem}
\label{r:Koszul}
    In the classical case of the Corona problem it is known that the Koszul complex provides a purely algebraic mechanism to arrive at resulting $\overline{\partial}$-equations.  This is well explained in \cite{MR2261424}*{Chapter 8, Appendix}.  Namely, when wanting to solve
    \begin{align*}
        \sum_{j=1}^n f_jg_j=1
    \end{align*}
    for holomorphic solutions $g_j$, one starts with solutions $\varphi_j$ that are smooth and bounded, and then perturbs these solutions to be holomorphic via:
    $$
    g_j=\varphi_j+\sum_{\overset{k=1}{k\neq j}}^n(b_{k,j}-b_{j,k})f_k
    $$
    where $\overline{\partial}b_{j,k}=\varphi_j\overline{\partial}\varphi_k$.  This last condition is what forces the holomorphicity of the solutions $g_j$.  
    
    It would be interesting to formulate the related Koszul complex in this setting to produce a purely algebraic mechanism to provide the holomorphic solutions $H_j$ and $K_j$ to \eqref{e:nNewBezout1} and \eqref{e:nNewBezout2}.
\end{rem}

\subsection{Wolff's Estimates and Bounded Solutions to \texorpdfstring{$\overline{\partial}$}{d bar}-Equations}
\label{s:Wolff}

To see that Wolff's Theorem applies, one simply verifies certain estimates on the $\overline{\partial}\beta_j$.  These are similar to the classical case and we provide details here to explain how they can be applied in the setting at hand.

The main tool to exploit is the following lemma of Wolff, see \cite{MR2261424}*{Chapter~8, Theorem~1.2}.  See also the monographs \cites{MR1669574,MR1419088} where this result is also provided.  We additionally utilize an observation of Treil to produce bounded solutions on the entire disk $\mathbb{D}$.

To state the result, recall that a measure $\mu$
 is a Carleson measure if there exists a constant $C$ such that
 $$
 \mu\left(Q_I\right)\leq C|I|\quad\forall I\subset \mathbb{T}, \quad Q_I=\{z=re^{i\theta}:e^{i\theta}\in I, 1-\vert I\vert\leq r<1\}.
 $$
 Call the best possible constant in the above the \textit{norm of the Carleson measure} $\mu$, written $\Vert \mu\Vert_{\mathcal{C}}$.  There is a close connection between Carleson measures and $BMO$ functions.  Namely, a $f\in BMOA(\mathbb{D})$ if and only if $d\mu_f(z):=\vert f'(z)\vert^2\log\frac{1}{\vert z\vert}dA(z)$ is a Carleson measure and moreover,
 $$\Vert \mu_f\Vert_{\mathcal{C}}\approx \Vert f\Vert_{BMOA}^2\lesssim\Vert f\Vert_{H^\infty(\mathbb{D})}^2.
 $$

 The result which then produces bounded solutions to equations of the form $\overline{\partial}b=G$ is given by the following:
\begin{thm}[Wolff, \cite{MR2261424}*{Chapter~8, Theorem~1.2}, and Treil, \cite{MR1183608}*{Theorem~3.1}]
\label{t:boundedsolutions}
    Suppose that $G(z)$ is bounded and $C^1(\mathbb{D})$.  Assume that $\vert G\vert^2\log\frac{1}{\vert z\vert}dA(z)$ and $\vert \partial G\vert \log\frac{1}{\vert z \vert}dA(z)$ are Carleson measures with norm $B_1$ and $B_2$ respectively, and that
    $$
    \left\vert\partial G\right\vert\leq\frac{B_3}{(1-\vert z\vert^2)^2}.
    $$  Then there is a function $b\in C(\overline{\mathbb{D}})\cap C^1(\mathbb{D})$ solving the equation $\overline{\partial} b=G$ with $\|b\|_{L^\infty(\mathbb{D})}\lesssim \sqrt{B_1}+B_2+B_3$.
\end{thm}

\subsubsection{Proof of Wolff's Estimates.}
The above result is applied when $G$ has the form \eqref{e:dbetarsj}, \eqref{e:dgammarsj}, \eqref{e:detarsj} or \eqref{e:detatildersj}.  The functions $G$ we are working with are of the form
\begin{itemize}
    \item $\overline{\partial}h_s\left(\hat{h}_rh_j+\hat{k}_rk_j\right)
        -\overline{\partial}h_j\left(\hat{h}_rh_s+\hat{k}_rk_s\right)
        -\overline{\partial}\hat{k}_r\left(h_sk_j-h_jk_s\right)$;
    \item $\overline{\partial}\hat{k}_s\left(\hat{h}_rh_j+\hat{k}_rk_j\right)
        -\overline{\partial}\hat{k}_r\left(\hat{h}_sh_j+\hat{k}_sk_j\right)
        -\overline{\partial}h_j\left(\hat{h}_r\hat{k}_s-\hat{h}_s\hat{k}_r\right)$;
    \item $- \overline{\partial}h_s\left(h_jk_r-h_rk_j\right)
        -\overline{\partial}h_r\left(h_sk_j-h_jk_s\right)
        -\overline{\partial}h_j\left(h_rk_s-h_sk_r\right)$;
    \item $- \overline{\partial}\hat{k}_s\left(\hat{h}_r\hat{k}_j-\hat{h}_j\hat{k}_r\right)
        -\overline{\partial}\hat{k}_r\left(\hat{h}_j\hat{k}_s-\hat{h}_s\hat{k}_j\right)
        -\overline{\partial}\hat{k}_j\left(\hat{h}_s\hat{k}_r-\hat{h}_r\hat{k}_s\right)$.
\end{itemize}
Since $\Vert h_j\Vert_{L^\infty(\mathbb{D})}$ and $\Vert k_j\Vert_{L^\infty(\mathbb{D})}$ are bounded by
$$
C_1(n,\delta):=\frac{1}{\delta^2}+\frac{n}{\delta^4},
$$ we see that it would ultimately suffice to have information about $\vert\overline{\partial} h_s\vert^2$, $\vert\overline{\partial} \hat{h}_s\vert^2$, $\vert\overline{\partial} k_s\vert^2$, $\vert\overline{\partial} \hat{k}_s\vert^2$, $1\leq s\leq n$, being controlled by $\displaystyle \sum_{j=1}^n\left(\vert F_j'\vert^2+\vert \hat{F}_j'\vert^2+\vert G_j'\vert^2+\vert \hat{F}_j'\vert^2\right)$.  Further for $G$ of the form above, we will have that $\vert \partial G\vert$ can be dominated by $\vert\overline{\partial}h_s\vert\vert\partial h_r\vert$, $\vert\overline{\partial}\hat{h}_s\vert\vert\partial h_r\vert$, $\vert\overline{\partial}h_s\vert\vert\partial \hat{h}_r\vert$, $\vert\overline{\partial}\hat{h}_s\vert\vert\partial \hat{h}_r\vert$, $\vert\overline{\partial}k_s\vert\vert\partial k_r\vert$, $\vert\overline{\partial}k_s\vert\vert\partial \hat{k}_r\vert$, $\vert\overline{\partial}\hat{k}_s\vert\vert\partial k_r\vert$, $\vert\overline{\partial}\hat{k}_s\vert\vert\partial \hat{k}_r\vert$, $\vert\overline{\partial}h_s\vert\vert\partial \hat{k}_r\vert$, $\vert\overline{\partial}\hat{h}_s\vert\vert\partial k_r\vert$, $\vert\overline{\partial}h_s\vert\vert\partial k_r\vert$ or $\vert\overline{\partial}\hat{h}_s\vert\vert\partial \hat{k}_r\vert$ and then the trivial inequality $xy\lesssim x^2+y^2$ will allow us to control $\vert\partial G\vert$ by
$$
\displaystyle \sum_{j=1}^n\left(\vert F_j'\vert^2+\vert \hat{F}_j'\vert^2+\vert G_j'\vert^2+\vert \hat{F}_j'\vert^2\right).
$$

This estimate would then give the claimed control on the Carleson measures and data in Theorem \ref{t:boundedsolutions} because $F_j,G_j\in H^\infty(\mathbb{D})$, the Carleson estimates follow from the containment $H^\infty(\mathbb{D})\subset BMOA(\mathbb{D})$ and the estimate

$$\vert f'(z)\vert\leq\frac{\Vert f\Vert_{H^\infty(\mathbb{D})}}{(1-\vert z\vert^2)}\ \ {\rm for\ all}\ \  f\in H^\infty(\mathbb{D}).
$$

 In particular, it will be the case that $B_1\lesssim C_1(n,\delta)^2 \left(\frac{n}{\delta^4}+\frac{n^2}{\delta^8}\right)$, $B_2\lesssim C_1(n,\delta) \left(\frac{n}{\delta^4}+\frac{n^2}{\delta^8}\right)$ and $B_3\lesssim C_1(n,\delta) \left(\frac{n}{\delta^4}+\frac{n^2}{\delta^8}\right)$.  These estimates are extremely coarse and one should be able to do much better by being more careful!

In terms of the estimates on $\vert\overline{\partial} h_s\vert^2$, $\vert\overline{\partial} \hat{h}_s\vert^2$, $\vert\overline{\partial} k_s\vert^2$, $\vert\overline{\partial} \hat{k}_s\vert^2$, $1\leq s\leq n$, from \eqref{e:hjfunctions} and \eqref{e:kjfunctions} we have
\begin{align*}
h_j & := \varphi_j+\sum_{\overset{r=1}{r\neq j}}^n\lambda_{r,j} F_r+\sum_{r=1}^n\mu_{r,j}G_r,\quad 1\leq j\leq n\\
k_j & := \psi_j+\sum_{\overset{r=1}{r\neq j}}^n\hat{\alpha}_{r,j}\hat{G}_r+\sum_{r=1}^n\hat{\mu}_{j,r}\hat{F}_r,\quad 1\leq j\leq n
\end{align*}
and so one needs resulting estimates on $\vert\overline{\partial}\varphi_j\vert$, $\vert\overline{\partial}\psi_j\vert$, $\vert\overline{\partial}\lambda_{r,j}\vert$, $\vert\overline{\partial}\alpha_{r,j}\vert$, $\vert\overline{\partial}\mu_{r,j}\vert$ and $\vert\overline{\partial}\hat{\mu}_{r,j}\vert$ for all $1\leq r,j\leq n$.  However, at this stage, the estimate are analogous to those contained in \cite{MR2261424}*{Chapter 8, Theorem 1.2} or \cite{MR1419088}*{Chapter 8, Section 2} and we sketch how to obtain the claimed control.

The estimate for $\vert\overline{\partial}\varphi_j\vert^2$ proceeds as follows.  We know that $\varphi_j(z)=\frac{\overline{F_j(z)}}{D(z)}$, which by the chain rule and simple estimates gives
\begin{align*}
\vert \overline{\partial}\varphi_j\vert^2 & = \frac{\vert D \overline{F_j'}-\overline{F_j}\,\overline{\partial}D\vert^2}{D^2}\\
& \lesssim \frac{n}{\delta^4}\sum_{j=1}^n \left(\vert F_j'\vert^2+\vert \hat{F}_j'\vert^2+\vert G_j'\vert^2+\vert \hat{G}_j'\vert^2\right)
\end{align*}
since
$\displaystyle\overline{\partial}D=\sum_{j=1}^nF_j \overline{F_j'}+G_j \overline{G_j'}$
can be estimated as
$$
\displaystyle\vert \overline{\partial} D\vert^2\leq n \sum_{j=1}^n \left(\vert F_j'\vert^2+\vert \hat{F}_j'\vert^2+\vert G_j'\vert^2+\vert \hat{G}_j'\vert^2\right).
$$
 The terms for $\lambda_{r,j}, \hat{\alpha}_{r,j}, \mu_{r,j}$ and $\hat{\mu}_{r,j}$ are done in a similar fashion but uses the function $\Delta(z)$ in \eqref{e:NewDenominator}.
 For these terms a similar estimate holds,
 but with
 $$\displaystyle\frac{n^{\frac{3}{2}}}{\delta^8}\sum_{j=1}^n \left(\vert F_j'\vert^2+\vert \hat{F}_j'\vert^2+\vert G_j'\vert^2+\vert \hat{G}_j'\vert^2\right),
  $$
  which when combined all together will give:
$$
\left\vert \overline{\partial} h_j\right\vert^2 \lesssim \left(\frac{n}{\delta^4}+\frac{n^2}{\delta^8}\right)\sum_{j=1}^n \left(\vert F_j'\vert^2+\vert \hat{F}_j'\vert^2+\vert G_j'\vert^2+\vert \hat{G}_j'\vert^2\right).
$$
The estimates for $\vert\partial h_s\vert^2$, $\vert\partial \hat{h}_s\vert^2$, $\vert\partial k_s\vert^2$, $\vert\partial \hat{k}_s\vert^2$, $1\leq s\leq n$ are handled analogously and details are omitted.

\section{Solution of the Slice Hyperholomorphic Corona Problem}
\label{s:QCorona}

In this section we show how Theorem ~\ref{t:NewCoronaS2} from Section~\ref{s:NewCorona} can be used to solve a resulting Corona problem in the slice hyperholomorphic setting.  Namely we will prove Theorem~\ref{t:QCorona}.  First we need to set up the appropriate notation and terminology.

\subsection{Preliminaries on Slice Hyperholomorphic Functions}
\label{s:SlicePrelims}

Much of the preliminary material recalled here can be found in \cites{MR2752913,MR3585395,MR3013643}.

We let $\mathbb{H}$ denote the set of quaternions.  This will be the set $\mathbb{R}^4$ endowed with a multiplicative structure.  More precisely, given $q=(x_0,x_1,x_2,x_3)\in\mathbb{R}^4$ we write:
$$
q=x_0e_0+x_1e_1+x_2e_2+x_3e_3
$$
where $\{e_0,e_1,e_2,e_3\}$ denotes the standard basis of $\mathbb{R}^4$, $e_0=1$ and the following multiplication rule for the other basis elements hold:
\begin{align*}
    e_1^2 &= e_2^2=e_3^2=-1\\
    e_1e_2 &= -e_2e_1=e_3\\
    e_2e_3 &= -e_3e_2=e_1\\
    e_3e_1 &= -e_1e_3=e_2.
\end{align*}
Multiplication between $q_1,q_2\in\mathbb{R}^4$ is then done in the obvious fashion subject to the rules given above.  The elements $\{e_1,e_2,e_3\}$ can be thought of as the imaginary units within this algebra.  Associated to $q\in\mathbb{H}$ define the conjugate as
$$
\overline{q}:=x_0e_0-x_1e_1-x_2e_2-x_3e_3
$$
and the norm of $q$ as $\vert q\vert^2:=\sqrt{q\overline{q}}=x_0^2+x_1^2+x_2^2+x_3^2.$  From this formula one can see that for $q\neq 0$ that the inverse to $q$ is given by $q^{-1}:=\vert q\vert^{-2} \overline{q}$.  We further define $\operatorname{Re}q=x_0$ and $\operatorname{Im}q=x_1e_1+x_2e_2+x_3e_3$.  For $q\in\mathbb{H}$ with $\operatorname{Im}q\neq 0$, one sees that $\frac{\operatorname{Im}q}{\vert \operatorname{Im}q\vert}$ is an imaginary unit.  From this, one sees that any $q\in\mathbb{H}$, with $q\notin\mathbb{R}$ can be written as $q=u+Iv$ for $u,v\in\mathbb{R}$. Let
$$
\mathbb{S}=\{q\in\mathbb{H}: q^2=-1\}
$$
denote the set of imaginary units.  For every $I\in\mathbb{S}$ set $\mathbb{C}_I$ the plane $\mathbb{R}+I\mathbb{R}$.  Finally, for a subset $\Omega\subset\mathbb{H}$, we will set $\Omega_I:=\Omega\cap \mathbb{C}_I$, the intersection of the complex plane $\mathbb{C}_I$ with the subset $\Omega$.

A domain $\Omega$ in $\mathbb{H}$ is called a \textit{slice domain} if for all $I\in\mathbb{S}$ the slice $\Omega_I$ is a domain of the complex plane $\mathbb{C}_I$.  The domain is called a \textit{axially symmetric domain} if for all $x,y\in\mathbb{R}$, $x+Iy\in\Omega$ implies $x+\mathbb{S}y\subset\Omega$.  The domain we will focus on in this section will be $\mathbb{B}$, the unit ball in $\mathbb{H}$;  namely, the set $\{q\in\mathbb{H}: \vert q\vert <1\}$.  This is an example of an axially symmetric slice domain in $\mathbb{H}$.

We will be interested in the hyperholomorphic (or regular) functions in the quaternionic setting.  A function $f:\Omega\subset\mathbb{H}$ is called \textit{slice hyperholomorphic} (or \textit{slice regular}) if for all $I\in\mathbb{S}$ the restriction $f_I:\Omega_I\to\mathbb{H}$ has continuous partial derivatives and is holomorphic on $\Omega_I$, namely
$$
\overline{\partial}_If_I(x+yI):=\frac{1}{2}\left(\frac{\partial}{\partial x}+I\frac{\partial}{\partial y}\right)f_I(x+yI)=0.
$$

A key result that relates slice hyperholomorphic functions to classical holomorphicity is contained in the ``Splitting Lemma.''
\begin{lm}[Splitting Lemma]
\label{l:splitlem}
If $f$ is a slice hyperholomorphic function on $\Omega$, then for every $I\in\mathbb{S}$ and every $J\in\mathbb{S}$ such that $J\perp I$, there exists two holomorphic functions $F,G:\Omega_I\to \mathbb{C}_I$ such that for every $z=x+Iy\in\Omega_I$ it holds that
$$
f_I(z)=F(z)+G(z)J.
$$
\end{lm}
The key point of this lemma is that allows the study of slice hyperholomorphicity by restricting to slices in the domain and then appealing to classical holomorphicity of a complex variable.

Interestingly, this can be reversed via the ``Extension Lemma'' which says the behavior on a slice can be used to go to the whole domain.
\begin{lm}[Extension Lemma]
\label{l:extlem}
Let $\Omega$ be an axially symmetric slice domain in $\mathbb{H}$ and choose $I\in\mathbb{S}$.  If $f_I:\Omega_I\to \mathbb{C}_I$ is holomorphic, then setting
\begin{equation}\label{RFormula}
f(q)=f(x+yL)=\frac{1}{2}\left(f_I(x+yI)+f_I(x-yI)\right)+L\frac{I}{2}\left(f_I(x-yI)-f_I(x+yI)\right)
\end{equation}
extends $f_I$ to a slice hyperholomorphic function $f:\Omega\to \mathbb{H}$.
\end{lm}
Note that any slice hyperholomorphic function $f$ on an axially symmetric slice domain $\Omega$ satisfies formula \eqref{RFormula} which is known as Representation Formula. 

There is a natural product on functions in this space, denoted by $\ast$.  This product is not commutative, but is associative and distributive.  It may be defined by:
$$
(f\ast g)(q):=f(q)g(f^{-1}(q)qf(q))
$$
when $f(q)\neq 0$ and zero when $f(q)=0$.

Let $f^{c}$ be the regular conjugate of $f$; the definition of $f^c$ can be given in general, but on a ball centered at the origin where $f(q)=\sum_{n=0}^\infty q^na_n$, $f^c$ is defined as $f^c(q)=\sum_{n=0}^\infty q^n\overline {a_n}$. Let $f^s(q):=(f^{c}\ast f)(q)=(f\ast f^{c})(q)$ denote the symmetrization of the function.

The $\ast$-inverse of a slice hyperholomorphic function $f$ is defined in the following way.  Define $T_f:\mathbb{B}\to\mathbb{B}$ by the formula $T_f(q):=f^{c}(q)^{-1}q f^{c}(q)$.  This preserves the modulus of the element $q$ and hence for any $q\in\mathbb{B}$ it will be true that $T_f(q)\in\mathbb{B}$.  Then $\cdot^{-\ast}$ is the operation given by
\begin{equation}
\label{e:*inverse}
f^{-\ast}(q)=f^{s}(q)^{-1} f^{c}(q).
\end{equation}
It is the case that $f^{-\ast}(q)$ is slice hyperholomorphic when it is defined.  It is known that both $f^{c}$ and $f^{s}$ are hyperholomorphic functions when $f$ is.

\subsection{Statement of \texorpdfstring{$H^\infty(\mathbb{B})$}{the Bounded Quaternionic Valued}-Corona Problem.}

The quaternionic Hardy spaces were introduced in \cite{MR3127378}, while a systematic treatment is in the book \cite{MR3585855}; further developments are in
\cite{MR3801294}. 
For $1\leq p< \infty$ we set:
$$
\left\Vert f\right\Vert_{H^{p}(\mathbb{B})}:=\sup_{I\in\mathbb{S}}\left\Vert f_I\right\Vert_{H^{p}(\mathbb{B}_I)}=\sup_{I\in\mathbb{S}}\lim_{r\to 1}\left(\frac{1}{2\pi}\int_{\mathbb{T}} \left\vert f(re^{I\theta}\right\vert^{p}d\theta\right)^{\frac{1}{p}},
$$
\noindent and for $p=\infty$, we have:
$$
\left\Vert f\right\Vert_{H^{\infty}(\mathbb{B})}=\sup_{q\in\mathbb{B}}\vert f(q)\vert=\sup_{I\in\mathbb{S}}\left\Vert f_I\right\Vert_{H^{\infty}(\mathbb{B}_I)}.
$$
We will not need the spaces $H^p(\mathbb{B})$ in this paper, but only the algebra $H^\infty(\mathbb{B})$.  It is immediate to observe that when we consider $\mathbb{B}\cap \mathbb{C}_I$ this is a disk $\mathbb{D}_I$ in the complex plane $\mathbb{C}_I$.

We make a couple of immediate remarks related to the boundedness properties.  First, one easily sees that if $f\in H^\infty(\mathbb{B})$ we have that $f^{c},f^{s}\in H^\infty(\mathbb{B})$, and the norm is related to the norm of $\Vert f\Vert_{H^\infty(\mathbb{B})}$.  More precisely, $\Vert f^{c}\Vert_{H^\infty(\mathbb{B})}=\Vert f\Vert_{H^\infty(\mathbb{B})}$ and $\Vert f^{s}\Vert_{H^\infty(\mathbb{B})}\leq \Vert f\Vert_{H^\infty(\mathbb{B})}^2$.  Second, by the Extension Lemma, Lemma \ref{l:extlem}, if we have $f_I$ bounded on a slice, then the extension will be bounded on $\mathbb{B}$.

Consider now $f_1,\ldots,f_n\in H^\infty(\mathbb{B})$.  We would like to study the Bezout equation for these functions, namely, we want to have functions $g_1,\ldots,g_n\in H^\infty(\mathbb{B})$ so that
\begin{equation}
\label{e:Bezout}
\tag{B}
\sum_{j=1}^n  (f_j \ast g_j)(q)=1.
\end{equation}
We focus on the situation of right multiplication by the $g_j$ because the case of left multiplication by the $g_j$, the quantitative condition of no common zeros is not sufficient as observed in \cite{Shelah_MSThesis}.  Indeed, if we consider the functions $f_1(q)=q-\frac{\imath}{2}$, $f_2(q)=\left(q-\frac{\imath}{2}\right)\jmath$ for $\imath,\jmath\in\mathbb{S}$ with $\imath\perp\jmath$.  Clearly we have $f_1\left(\frac{\imath}{2}\right)=f_2\left(\frac{\imath}{2}\right)=0$.  Yet if we take $g_1(q)=\imath$ and $g_2(q)=\imath\jmath$, then
$$
g_1(q)\ast f_1(q)+g_2(q)\ast f_2(q)=q\left(\imath+\imath\jmath^2\right)-\frac{\imath^2+(\imath\jmath)^2}{2}=1.
$$
When working with right multiplication by the $g_j$ this challenge goes away.

Observe that \eqref{e:Bezout} implies that $f_j$  have no common zeros in a quantitative sense.  Indeed, we have:
\begin{align*}
    1 & = \left\vert \sum_{j=1}^n f_j\ast g_j(q)\right\vert \leq \sum_{j=1}^n\left\vert f_j\ast g_j(q)\right\vert = \sum_{j=1}^n \left\vert f_j(q)\right\vert\left\vert g_j(f_j(q)^{-1}qf_j(q))\right\vert\\
    &\leq \left(\sum_{j=1}^n\left\vert f_j(q)\right\vert^2\right)^{\frac{1}{2}}\left(\sum_{j=1}^n\left\vert g_j(f_j(q)^{-1}qf_j(q))\right\vert^2\right)^{\frac{1}{2}} \leq \sqrt{n}\max_{j}\Vert g_j\Vert_{H^\infty(\mathbb{B})}\left(\sum_{j=1}^n\left\vert f_j(q)\right\vert^2\right)^{\frac{1}{2}}.
\end{align*}
The computation above proves that a necessary condition for \eqref{e:Bezout} is the following: there exists $\delta>0$ such that:
\begin{equation}
    \label{e:Carl}
    \tag{C}
    0<\delta^2\leq \sum_{j=1}^{n} \left\vert f_j(q)\right\vert^2\quad\forall q\in\mathbb{B}.
\end{equation}
The Corona problem in the setting of $H^\infty(\mathbb{B})$ is to prove that \eqref{e:Carl} implies \eqref{e:Bezout}.  We show this now in the following subsections.

\subsection{The Case of One Generator}
Here we argue why Theorem~\ref{t:QCorona} is true in the case of one function.  We give. direct proof as it is also possible to deduce this from Thereom \ref{t:NewCoronaS2}.  This will follow by elementary properties of slice hyperholomorphic functions as given above.

Suppose we have a function $f\in H^\infty(\mathbb{B})$ such that $\vert f(q)\vert^2\geq\delta^2>0$ for all $q\in\mathbb{B}$.  We want to show that there is a function $g\in H^\infty(\mathbb{B})$ so that $f\ast g (q)=1$ for all $q\in\mathbb{B}$ and
$$
\Vert g\Vert_{H^{\infty}(\mathbb{B})}\leq C(\delta)\Vert f\Vert_{H^{\infty}(\mathbb{B})}.
$$

To do this, we will use the $\ast$-inverse. Observe that
$$T_f(q):=f^{c}(q)^{-1}q f^{c}(q)$$ is well defined for all $q\in\mathbb{B}$ since $\vert f(q)\vert^2\geq\delta^2>0$ for all $q\in\mathbb{B}$.  Then define
$$
f^{-\ast}(q)=f^{s}(q)^{-1} f^{c}(q).
$$
It is the case that $f^{-\ast}(q)$ is hyperholomorphic and we observe that
$$
f\ast f^{-\ast}(q)=1\quad\forall q\in\mathbb{B},
$$
and so we have solved the Bezout equation pointwise for all $q\in\mathbb{B}$.  Note this argument is essentially \cite{MR3013643}*{Theorem~5.3}; see also \cite{MR3801294}*{Proposition~2.15}.

Next, we observe that $f^{-\ast}\in H^\infty(\mathbb{B})$ when $f\in H^\infty(\mathbb{B})$ and $1\geq \vert f(q)\vert^2\geq\delta^2>0$ for all $q\in\mathbb{B}$.  Indeed, since by the hypothesis that $1\geq \vert f(q)\vert^2\geq\delta>0$ for all $q\in\mathbb{B}$, we have that the same holds true for $\vert f^c(q)\vert^2\geq\delta^2>0$ for all $q\in\mathbb{B}$, in fact writing $q=x+Iy$, by the Representation Formula we have $f(x+Iy)=a+Ib$ and $|f^c(x+Iy)|=|\bar a+I\bar b|=|\overline{a-bI}|=a-bI=a+Jb=f(x+Jy)$ for some $J\in\mathbb S$.
By the definition of $f^{s}=f\ast f^{c}=f^{c}\ast f$ it is the case that $\vert f^{s}(q)\vert^2\geq\delta^4>0$ for all $q\in\mathbb{B}$ too.  Then, with this fact, we have
$$\Vert f^{-\star}\Vert_{H^\infty(\mathbb{B})}\leq\frac{1}{\delta^2} \Vert f\Vert_{H^\infty(\mathbb{B})}.
$$

Hence, for any $f\in H^\infty(\mathbb{B})$ we have a function $g\in H^{\infty}(\mathbb{B})$ so that $f\ast g(q)=1$ for all $q\in\mathbb{B}$, and $\Vert g\Vert_{H^{\infty}(\mathbb{B})}\leq C(\delta)\Vert f\Vert_{H^{\infty}(\mathbb{B})}$ with  $C(\delta)=\frac{1}{\delta^2}$.

\subsection{Proof of The Quaternionic Corona Theorem: Theorem~\ref{t:QCorona}}

Here we want to prove the following main result of this paper, restated from the Introduction:
\begin{thm}
    Let $0<\delta<1$ and $n\in\mathbb{N}$.  Suppose that $f_1,\ldots, f_n\in H^\infty(\mathbb{B})$ are slice hyperholomorphic functions such that
    $$
        0<\delta^2\leq \sum_{j=1}^n\left\vert f_j(q)\right\vert^2\leq 1\quad\forall q\in\mathbb{B}.
    $$
Then there are $g_1,\ldots, g_n$ such that
    $$
    \sum_{j=1}^n (f_j \ast g_j)(q)=1\quad\forall q\in\mathbb{B}
    $$
    and
    $$
    \displaystyle\max_{1\leq j\leq n} \Vert g_j\Vert_{H^\infty(\mathbb{B})}\leq C(\delta,n):= \frac{1}{\delta^2}+\frac{n}{\delta^4}+\left(\frac{1}{\delta^2}
    +\frac{n}{\delta^4}\right)\left(\frac{n}{\delta^4}+\frac{n^2}{\delta^8}\right).
    $$
\end{thm}

We follow the initial ideas and reductions of Gentili, Sarfatti and Struppa in \cite{MR3621101} where they studied ideals of slice hyperholomorphic functions in a quaternionic variable.  Their idea was to utilize the Splitting Lemma to pass to a system of Bezout-type equations and then solve that system of equations locally and use sheaf theory to patch together local solutions to have global solutions.  This method does not allow to have bounded information, yet nevertheless the approach they essentially follow is important for the approach we take.  We recall the setup and formulation here now.

Suppose now that we have $f_1,\ldots,f_n\in H^\infty(\mathbb{B})$ such that
$$
0<\delta^2\leq\sum_{j=1}^n \vert{f_j(q)}\vert^2\leq 1\quad\forall q\in \mathbb{B},
$$
recall that this is condition \eqref{e:Carl}.

We want to solve for functions $g_j\in H^\infty(\mathbb{B})$ so that
\begin{equation}
\label{e:temp}
\sum_{j=1}^{n} (f_j\ast g_j)(q)=1\quad\forall q\in\mathbb{B}.
\end{equation}
Via the Splitting Lemma, Lemma \ref{l:splitlem}, for any $J\perp I$ we can write for $1\leq j\leq n$:
\begin{eqnarray*}
    f_j(x+Iy) & = & F_j(x,y)+G_j(x,y)J\\
    g_j(x+Iy) & = & H_j(x,y)+K_j(x,y)J.
\end{eqnarray*}
In the slice $z=x+Iy\in\mathbb{D}_I$ we have that
$$
f_j(z)\ast g_j(z)=F_j(z)H_j(z)-G_j(z)\hat{K}_j(z)+\left(F_j(z)K_j(z)+G_j(z)\hat{H}_j(z)\right)J.
$$
Collecting and equating real and imaginary  parts (with respect to the complex unit $J$) of \eqref{e:temp} we arrive at the system of equations:
\begin{align}
    \sum_{j=1}^n F_jH_j-\sum_{j=1}^nG_j\hat{K}_j & =1\label{e:tempBezout1}\\
    \sum_{j=1}^n F_jK_j+\sum_{j=1}^nG_j\hat{H}_j & =0\label{e:tempBezout2}
\end{align}
that must be solved on the domain $\mathbb{D}_I$ (which is the unit disk in the complex plane $\mathbb{C}_I$).  The system of equations \eqref{e:tempBezout1} and \eqref{e:tempBezout2} is exactly \eqref{e:nNewBezout1} and \eqref{e:nNewBezout2}.  One observes that condition \eqref{e:Carl} implies a similar quantitative condition on the zeros of $F_j$ and $G_j$, which is the same as \eqref{e:nCarl}.

The system is exactly the same system as studied in Theorem~\ref{t:NewCorona}, and since we have that the data $F_j$ and $G_j$ satisfy \eqref{e:nCarl}, there exists solutions $H_j$ and $K_j$ in $H^{\infty}(\mathbb{D}_I)$.  Hence we have bounded holomorphic functions $H_j$ and $K_j$ on the whole disk $\mathbb{D}_I$ (the unit disk of the slice we are working on) that satisfy \eqref{e:nNewBezout1} and \eqref{e:nNewBezout2} and hence to \eqref{e:temp}.  By the Extension Lemma, Lemma~\ref{l:extlem}, we have a solution to \eqref{e:temp} that is bounded on all of $\mathbb{B}$.  By the estimates from Theorem~\ref{t:NewCorona} we have that $\Vert g_j\Vert_{H^\infty(\mathbb{B})}\leq C(\delta, n)$ appearing in Theorem \ref{t:NewCorona}.  This argument completes the proof of Theorem \ref{t:QCorona}.

\section{Concluding Remarks}
\label{s:Conclude}

There are many possible extensions of the results in this paper, especially Theorem \ref{t:NewCorona}.  This is because the Corona result of Carleson from \cite{MR0141789} has many possible connections to function theory and operator theory.  Additionally, the results from Theorem \ref{t:QCorona} would have related extensions as outlined below.

We point to possible lines of inquiry to extend the results in this paper:
\begin{enumerate}
    \item In the case of the classical Corona theorem, it is known that one can take the constant $C(\delta,n)$ independent of $n$.  In particular, this has an interpretation of an infinite number of generators $f_j\in H^\infty(\mathbb{D})$ having no common zeros leading to solutions $g_j\in H^\infty(\mathbb{D})$ so that $\displaystyle \sum_{j=1}^\infty f_jg_j=1$.  See for example \cites{MR0222701, MR0570865, MR0629839, Uchiyama, MR0595742} that explore aspects of this question in the classical case.  A natural question is the same for Theorems \ref{t:NewCorona} and \ref{t:QCorona}.  Are these results true when $n=\infty?$  There are additional questions about analogs of Theorems \ref{t:NewCorona} and \ref{t:QCorona} when the resulting Corona data is matrix or operator valued.  See \cites{MR0981054,MR0222701,MR2449054, MR2106344} for possible motivation as to questions one could explore in this direction.

    Additionally, in this direction one can seek to determine better information about the constant $C(\delta,n)$ and the absolute constants there.  See for example \cites{MR2158178, MR1887635, MR2213732, MR2317961, MR0990193, MR0595742, MR2105955} where estimates on the absolute constants are explored in the classical case.

    \item Again in the classical case it is known that solving a version of the Corona problem with a Hilbert space target and some related norm control, then this would be sufficient to solve the problem with bounded analytic data.  This is related to the so called Toeplitz Corona Theorem and the complete Pick property of the Hardy space $H^2(\mathbb{D})$.  What is the analog of this question for the space $H^\infty(\mathbb{B})$?  Is there a Hilbert space version of Theorem \ref{t:NewCorona}?  See the results \cites{MR0383098, MR0482355, MR1882259, MR3857903} where questions of this type are studied for $H^\infty(\mathbb{D})$ and $H^2(\mathbb{D})$.

    \item Corona theorems and results such as Theorems \ref{t:NewCorona} and \ref{t:QCorona} have an algebraic interpretation as an ideal generated by a collection of functions generating the whole algebra.  One can then generalize this question to study other ideal membership questions and seek necessary and sufficient conditions so that a function $h$ belongs to an ideal generated by some collection of functions $f_1,\ldots, f_n$.  This has been well studied in the classical case, see for example the following, \cites{MR0595742,MR3283648,MR1297916,MR0812323, MR1033435, MR2362422, MR2085166, MR3857903}.  What are analogs of Theorems \ref{t:NewCorona} and \ref{t:QCorona} of these ideal questions?

    \item The Corona problem for multiplier algebras of analytic functions on spaces of holomorphic functions on the disk has also been well studied; see for example any of the following papers \cites{MR3077143,MR2057771,MR1651405}.  Equations \eqref{e:nNewBezout1} and \eqref{e:nNewBezout2} of course make sense in these multiplier algebras and the general computations contained in this paper should be extendable to those settings.  As a concrete question can one prove that for the multiplier algebra of the Dirichlet space has an analog of Theorem \ref{t:NewCorona}?

    \item The algebra $H^\infty(\mathbb{D})$ has been well studied and many other algebraic properties have been studied.  In particular properties like the Bass stable rank and topological stable rank have been studied in \cite{MR1183608} and \cite{MR1302836}.  Are similar results true for $H^\infty(\mathbb{B})$?  What are the interpretations of the stable rank conditions for the equations \eqref{e:nNewBezout1} and \eqref{e:nNewBezout2}?

\end{enumerate}


\begin{bibdiv}
\begin{biblist}

\bib{MR1882259}{book}{
   author={Agler, Jim},
   author={McCarthy, John E.},
   title={Pick interpolation and Hilbert function spaces},
   series={Graduate Studies in Mathematics},
   volume={44},
   publisher={American Mathematical Society, Providence, RI},
   date={2002},
   pages={xx+308}
}

\bib{MR3585855}{book}{
    AUTHOR = {Daniel Alpay},
     AUTHOR = {Fabrizio Colombo},
      AUTHOR = {Irene Sabadini},
     TITLE = {Slice hyperholomorphic {S}chur analysis},
    SERIES = {Operator Theory: Advances and Applications},
    VOLUME = {256},
 PUBLISHER = {Birkh\"{a}user/Springer, Cham},
      YEAR = {2016},
     PAGES = {xii+362},
      ISBN = {978-3-319-42513-9; 978-3-319-42514-6},
       DOI = {10.1007/978-3-319-42514-6},
       URL = {https://doi.org/10.1007/978-3-319-42514-6},
}
	
\bib{MR3127378}{article}{
    AUTHOR = {Daniel Alpay},
      AUTHOR = {Fabrizio Colombo},
       AUTHOR = {Irene Sabadini},
     TITLE = {Pontryagin-de {B}ranges-{R}ovnyak spaces of slice
              hyperholomorphic functions},
   JOURNAL = {J. Anal. Math.},
    VOLUME = {121},
      YEAR = {2013},
     PAGES = {87--125},
      ISSN = {0021-7670},
       URL = {https://doi-org.chapman.idm.oclc.org/10.1007/s11854-013-0028-8},
}

\bib{MR1419088}{book}{
   author={Mats Andersson},
   title={Topics in complex analysis},
   series={Universitext},
   publisher={Springer-Verlag, New York},
   date={1997},
   pages={viii+157},
   isbn={0-387-94754-X}
}

\bib{MR0990193}{article}{
   author={Mats Andersson},
   title={The corona theorem for matrices},
   journal={Math. Z.},
   volume={201},
   date={1989},
   number={1},
   pages={121--130}
}

\bib{MR0383098}{article}{
   author={William Arveson},
   title={Interpolation problems in nest algebras},
   journal={J. Functional Analysis},
   volume={20},
   date={1975},
   number={3},
   pages={208--233}
}

\bib{MR0812323}{article}{
   author={Jean Bourgain},
   title={On finitely generated closed ideals in $H^\infty(D)$},
   journal={Ann. Inst. Fourier (Grenoble)},
   volume={35},
   date={1985},
   number={4},
   pages={163--174}
}

\bib{MR0141789}{article}{
   author={Carleson, Lennart},
   title={Interpolations by bounded analytic functions and the corona
   problem},
   journal={Ann. of Math. (2)},
   volume={76},
   date={1962},
   pages={547--559}
}

\bib{MR1297916}{article}{
   author={Cegrell, Urban},
   title={Generalisations of the corona theorem in the unit disc},
   journal={Proc. Roy. Irish Acad. Sect. A},
   volume={94},
   date={1994},
   number={1},
   pages={25--30}
}

\bib{MR1033435}{article}{
   author={Cegrell, Urban},
   title={A generalization of the corona theorem in the unit disc},
   journal={Math. Z.},
   volume={203},
   date={1990},
   number={2},
   pages={255--261}
}


	

\bib{MR3967697}{book}{
    AUTHOR = {Fabrizio Colombo},
          AUTHOR = {Jonathan Gantner},
     TITLE = {Quaternionic closed operators, fractional powers and
              fractional diffusion processes},
    SERIES = {Operator Theory: Advances and Applications},
    VOLUME = {274},
 PUBLISHER = {Birkh\"{a}user/Springer, Cham},
      YEAR = {2019},
     PAGES = {viii+322},
      ISBN = {978-3-030-16408-9; 978-3-030-16409-6},
       DOI = {10.1007/978-3-030-16409-6},
       URL = {https://doi.org/10.1007/978-3-030-16409-6},
}
		





\bib{MR3887616}{book}{
 AUTHOR = {Fabrizio Colombo},
          AUTHOR = {Jonathan Gantner},
    AUTHOR = {David P. Kimsey},
     TITLE = {Spectral theory on the {S}-spectrum for quaternionic
              operators},
    SERIES = {Operator Theory: Advances and Applications},
    VOLUME = {270},
 PUBLISHER = {Birkh\"{a}user/Springer, Cham},
      YEAR = {2018},
     PAGES = {ix+356},
      ISBN = {978-3-030-03073-5; 978-3-030-03074-2},
}



	

\bib{MR4502763}{article}{
    AUTHOR = {Fabrizio Colombo},
       AUTHOR = {Antonino De Martino},
       AUTHOR = {Stefano Pinton},
       AUTHOR = {Irene Sabadini},
     TITLE = {Axially harmonic functions and the harmonic functional
              calculus on the {$S$}-spectrum},
   JOURNAL = {J. Geom. Anal.},
    VOLUME = {33},
      YEAR = {2023},
    NUMBER = {1},
     PAGES = {Paper No. 2, 54},
      ISSN = {1050-6926},
       URL = {https://doi-org.chapman.idm.oclc.org/10.1007/s12220-022-01062-3},
}


\bib{MR2555912}{article}{
   author={Colombo, Fabrizio},
   author={Gentili, Graziano},
   author={Sabadini, Irene},
   author={Struppa, Daniele},
   title={Extension results for slice regular functions of a quaternionic
   variable},
   journal={Adv. Math.},
   volume={222},
   date={2009},
   number={5},
   pages={1793--1808}
}

\bib{MR2089988}{book}{
    AUTHOR = {Fabrizio Colombo},
    AUTHOR = {Irene Sabadini},
   AUTHOR = {Franciscus Sommen},
     AUTHOR = {Daniele C. Struppa},
     TITLE = {Analysis of {D}irac systems and computational algebra},
    SERIES = {Progress in Mathematical Physics},
    VOLUME = {39},
 PUBLISHER = {Birkh\"{a}user Boston, Inc., Boston, MA},
      YEAR = {2004},
     PAGES = {xiv+332},
      ISBN = {0-8176-4255-2},
       DOI = {10.1007/978-0-8176-8166-1},
       URL = {https://doi.org/10.1007/978-0-8176-8166-1},
}
	
\bib{MR2752913}{book}{
    AUTHOR = {Colombo, Fabrizio},
    author={Sabadini, Irene},
   author={Struppa, Daniele C.},
     TITLE = {Noncommutative functional calculus},
    SERIES = {Progress in Mathematics},
    VOLUME = {289},
      NOTE = {Theory and applications of slice hyperholomorphic functions},
 PUBLISHER = {Birkh\"{a}user/Springer Basel AG, Basel},
      YEAR = {2011},
     PAGES = {vi+221},
      ISBN = {978-3-0348-0109-6},
       DOI = {10.1007/978-3-0348-0110-2},
       URL = {https://doi.org/10.1007/978-3-0348-0110-2},
}


\bib{MR3585395}{book}{
    AUTHOR = {Colombo, Fabrizio},
    author={Sabadini, Irene},
   author={Struppa, Daniele C.},
     TITLE = {Entire slice regular functions},
    SERIES = {SpringerBriefs in Mathematics},
 PUBLISHER = {Springer, Cham},
      YEAR = {2016},
     PAGES = {v+118},
      ISBN = {978-3-319-49264-3; 978-3-319-49265-0},
       DOI = {10.1007/978-3-319-49265-0},
       URL = {https://doi.org/10.1007/978-3-319-49265-0},
}

\bib{MR4240465}{book}{
    AUTHOR = {Colombo, Fabrizio},
    author={Sabadini, Irene},
    author={Struppa, Daniele C.},
     TITLE = {Michele {S}ce's works in hypercomplex analysis-a translation
              with commentaries},
 PUBLISHER = {Birkh\"{a}user/Springer, Cham},
      YEAR = {[2020] \copyright 2020},
     PAGES = {vi+122},
      ISBN = {978-3-030-50215-7; 978-3-030-50216-4},
       DOI = {10.1007/978-3-030-50216-4},
       URL = {https://doi.org/10.1007/978-3-030-50216-4},
}

\bib{MR2684426}{article}{
   author={Colombo, Fabrizio},
   author={Sabadini, Irene},
   author={Struppa, Daniele C.},
   title={An extension theorem for slice monogenic functions and some of its
   consequences},
   journal={Israel J. Math.},
   volume={177},
   date={2010},
   pages={369--389},
}

\bib{MR2520116}{article}{
   author={Colombo, Fabrizio},
   author={Sabadini, Irene},
   author={Struppa, Daniele C.},
   title={Slice monogenic functions},
   journal={Israel J. Math.},
   volume={171},
   date={2009},
   pages={385--403},
}


\bib{MR3077143}{article}{
   author={Costea, \c Serban},
   author={Sawyer, Eric T.},
   author={Wick, Brett D.},
   title={The corona theorem for the Drury-Arveson Hardy space and other
   holomorphic Besov-Sobolev spaces on the unit ball in $\mathbb{C}^n$},
   journal={Anal. PDE},
   volume={4},
   date={2011},
   number={4},
   pages={499--550}
}

\bib{MR3801294}{article}{
   author={de Fabritiis, Chiara},
   author={Gentili, Graziano},
   author={Sarfatti, Giulia},
   title={Quaternionic Hardy spaces},
   journal={Ann. Sc. Norm. Super. Pisa Cl. Sci. (5)},
   volume={18},
   date={2018},
   number={2},
   pages={697--733},
   issn={0391-173X}
}


\bib{MR4674273}{article}{
    AUTHOR = {Antonino De Martino},
       AUTHOR = {Stefano Pinton},
     TITLE = {Properties of a polyanalytic functional calculus on the
              {$S$}-spectrum},
   JOURNAL = {Math. Nachr.},
    VOLUME = {296},
      YEAR = {2023},
    NUMBER = {11},
     PAGES = {5190--5226},
      ISSN = {0025-584X},
       URL = {https://doi-org.chapman.idm.oclc.org/10.1002/mana.202200318},
}



\bib{MR4576314}{article}{
   AUTHOR = {Antonino De Martino},
       AUTHOR = {Stefano Pinton},
     TITLE = {A polyanalytic functional calculus of order 2 on the
              {$S$}-spectrum},
   JOURNAL = {Proc. Amer. Math. Soc.},
    VOLUME = {151},
      YEAR = {2023},
    NUMBER = {6},
     PAGES = {2471--2488},
      ISSN = {0002-9939},
       URL = {https://doi-org.chapman.idm.oclc.org/10.1090/proc/16285},
}	


\bib{MR4741057}{article}{
AUTHOR = {Antonino De Martino},
       AUTHOR = {Stefano Pinton},
    AUTHOR = {Peter Schlosser},
     TITLE = {The harmonic {$H^\infty$}-functional calculus based on the
              {$S$}-spectrum},
   JOURNAL = {J. Spectr. Theory},
    VOLUME = {14},
      YEAR = {2024},
    NUMBER = {1},
     PAGES = {121--162},
      ISSN = {1664-039X},
       URL = {https://doi-org.chapman.idm.oclc.org/10.4171/jst/492},
}

\bib{MR1169463}{book}{
    AUTHOR = {Richard Delanghe},
       AUTHOR = {Frank Sommen},
        AUTHOR = {Vladimir and Sou\v{c}ek},
     TITLE = {Clifford algebra and spinor-valued functions},
    SERIES = {Mathematics and its Applications},
    VOLUME = {53},
      NOTE = {A function theory for the Dirac operator,
              Related REDUCE software by F. Brackx and D. Constales,
              With 1 IBM-PC floppy disk (3.5 inch)},
 PUBLISHER = {Kluwer Academic Publishers Group, Dordrecht},
      YEAR = {1992},
     PAGES = {xviii+485},
      ISBN = {0-7923-0229-X},
       DOI = {10.1007/978-94-011-2922-0},
       URL = {https://doi.org/10.1007/978-94-011-2922-0},
}


\bib{MR3329539}{article}{
   author={Douglas, Ronald G.},
   author={Krantz, Steven G.},
   author={Sawyer, Eric T.},
   author={Treil, Sergei},
   author={Wick, Brett D.},
   title={A history of the corona problem},
   conference={
      title={The corona problem},
   },
   book={
      series={Fields Inst. Commun.},
      volume={72},
      publisher={Springer, New York},
   },
   isbn={978-1-4939-1254-4},
   isbn={978-1-4939-1255-1},
   date={2014},
   pages={1--29}
}



\bib{MR0222701}{article}{
   author={Fuhrmann, Paul A.},
   title={On the corona theorem and its application to spectral problems in
   Hilbert space},
   journal={Trans. Amer. Math. Soc.},
   volume={132},
   date={1968},
   pages={55--66}
}

\bib{MR2261424}{book}{
   author={Garnett, John B.},
   title={Bounded analytic functions},
   series={Graduate Texts in Mathematics},
   volume={236},
   edition={1},
   publisher={Springer, New York},
   date={2007},
   pages={xiv+459},
   isbn={978-0-387-33621-3},
   isbn={0-387-33621-4}
}

\bib{MR3621101}{article}{
   author={Gentili, Graziano},
   author={Sarfatti, Giulia},
   author={Struppa, Daniele C.},
   title={Ideals of regular functions of a quaternionic variable},
   journal={Math. Res. Lett.},
   volume={23},
   date={2016},
   number={6},
   pages={1645--1663}
}

\bib{MR3013643}{book}{
   author={Gentili, Graziano},
   author={Stoppato, Caterina},
   author={Struppa, Daniele C.},
   title={Regular functions of a quaternionic variable},
   series={Springer Monographs in Mathematics},
   publisher={Springer, Heidelberg},
   date={2013},
   pages={x+185},
   isbn={978-3-642-33870-0},
   isbn={978-3-642-33871-7},
   doi={10.1007/978-3-642-33871-7},
}

\bib{MR2353257}{article}{
   author={Gentili, Graziano},
   author={Struppa, Daniele C.},
   title={A new theory of regular functions of a quaternionic variable},
   journal={Adv. Math.},
   volume={216},
   date={2007},
   number={1},
   pages={279--301}
}

\bib{MR3496962}{book}{
AUTHOR = {Klaus G\"{u}rlebeck},
AUTHOR = {Klaus Habetha},
AUTHOR = {Wolfgang Spr\"{o}\ss ig},
     TITLE = {Application of holomorphic functions in two and higher
              dimensions},
 PUBLISHER = {Birkh\"{a}user/Springer, [Cham]},
      YEAR = {2016},
     PAGES = {xv+390},
      ISBN = {978-3-0348-0962-7; 978-3-0348-0964-1},
       DOI = {10.1007/978-3-0348-0964-1},
       URL = {https://doi.org/10.1007/978-3-0348-0964-1},
}

\bib{MR3857903}{article}{
   author={Hartz, Michael},
   author={Wick, Brett D.},
   title={Ideal membership in $H^{\infty}$: Toeplitz corona approach},
   journal={Integral Equations Operator Theory},
   volume={90},
   date={2018},
   number={6},
   pages={Paper No. 66, 16}
}

\bib{MR3283648}{article}{
   author={Holloway, Caleb D.},
   author={Trent, Tavan T.},
   title={Wolff's theorem on ideals for matrices},
   journal={Proc. Amer. Math. Soc.},
   volume={143},
   date={2015},
   number={2},
   pages={611--620}
}

\bib{MR2069293}{book}{
    AUTHOR = {Brian Jefferies},
     TITLE = {Spectral properties of noncommuting operators},
    SERIES = {Lecture Notes in Mathematics},
    VOLUME = {1843},
 PUBLISHER = {Springer-Verlag, Berlin},
      YEAR = {2004},
     PAGES = {viii+184},
      ISBN = {3-540-21923-4},
}


\bib{MR1669574}{book}{
   author={Koosis, Paul},
   title={Introduction to $H_p$ spaces},
   series={Cambridge Tracts in Mathematics},
   volume={115},
   edition={2},
   note={With two appendices by V. P. Havin [Viktor Petrovich Khavin]},
   publisher={Cambridge University Press, Cambridge},
   date={1998},
   pages={xiv+289}
}

\bib{MR0827223}{book}{
   author={Nikol\cprime ski\u i, Nikolai K.},
   title={Treatise on the shift operator},
   series={Grundlehren der mathematischen Wissenschaften [Fundamental
   Principles of Mathematical Sciences]},
   volume={273},
   note={Spectral function theory;
   With an appendix by S. V. Hru\v s\v cev [S. V. Khrushch\"ev] and V. V.
   Peller;
   Translated from the Russian by Jaak Peetre},
   publisher={Springer-Verlag, Berlin},
   date={1986},
   pages={xii+491}
}

\bib{MR2085166}{article}{
   author={Pau, Jordi},
   title={On a generalized corona problem on the unit disc},
   journal={Proc. Amer. Math. Soc.},
   volume={133},
   date={2005},
   number={1},
   pages={167--174}
}
\bib{MR3930595}{book}{
    AUTHOR = {Tao Qian},
      AUTHOR = {Pengtao Li},
     TITLE = {Singular integrals and {F}ourier theory on {L}ipschitz
              boundaries},
 PUBLISHER = {Science Press Beijing, Beijing; Springer, Singapore},
      YEAR = {2019},
     PAGES = {xv+306},
      ISBN = {978-981-13-6499-0; 978-981-13-6500-3},
       DOI = {10.1007/978-981-13-6500-3},
       URL = {https://doi.org/10.1007/978-981-13-6500-3},
}



\bib{MR0210910}{article}{
   author={Rajeswara Rao, K. V.},
   title={On a generalized corona problem},
   journal={J. Analyse Math.},
   volume={18},
   date={1967},
   pages={277--278}
}

\bib{MR0570865}{article}{
   author={Rosenblum, Marvin},
   title={A corona theorem for countably many functions},
   journal={Integral Equations Operator Theory},
   volume={3},
   date={1980},
   number={1},
   pages={125--137}
}

\bib{Saracco}{article}{
    author={Saracco, Alberto},
    publisher={John Wiley & Sons, Ltd},
    isbn={9781119414421},
    title={The Corona Problem, Carleson Measures, and Applications},
    booktitle={Mathematical Analysis and Applications},
    pages={709-730},
    year={2018}
}

\bib{MR0482355}{article}{
   author={ C. F. Schubert},
   title={The corona theorem as an operator theorem},
   journal={Proc. Amer. Math. Soc.},
   volume={69},
   date={1978},
   number={1},
   pages={73--76}
}


\bib{Shelah_MSThesis}{article}{
   author={Yonatan Shelah},
   title={Quaternionic Wiener algebras, factorization and applications},
   journal={MS Thesis, Tel Aviv University},
   date={2016},
   pages={1-32}
}

\bib{MR0916722}{article}{
   author={Sibony, Nessim},
   title={Probl\`eme de la couronne pour des domaines pseudoconvexes \`a{}
   bord lisse},
   language={French},
   journal={Ann. of Math. (2)},
   volume={126},
   date={1987},
   number={3},
   pages={675--682}
}

\bib{MR1302836}{article}{
   author={Su\'arez, Fernando Daniel},
   title={Cohomological stabilization of maximal ideal spaces in Banach
   algebras},
   journal={J. Algebra},
   volume={170},
   date={1994},
   number={1},
   pages={165--183}
}



\bib{MR0629839}{article}{
   author={Tolokonnikov, Vadim. A.},
   title={Estimates in the Carleson corona theorem, ideals of the algebra
   $H\sp{\infty }$, a problem of Sz.-Nagy},
   note={Investigations on linear operators and the theory of functions, XI},
   language={Russian, with English summary},
   journal={Zap. Nauchn. Sem. Leningrad. Otdel. Mat. Inst. Steklov. (LOMI)},
   volume={113},
   date={1981},
   pages={178--198, 267}
}


\bib{MR0595742}{article}{
   author={Tolokonnikov, Vadim. A.},
   title={Estimates in Carleson's corona theorem and finitely generated
   ideals of the algebra $H\sp{\infty }$},
   language={Russian},
   journal={Funktsional. Anal. i Prilozhen.},
   volume={14},
   date={1980},
   number={4},
   pages={85--86}
}



\bib{MR2362422}{article}{
   author={Treil, Sergei},
   title={The problem of ideals of $H^\infty$: beyond the exponent $3/2$},
   journal={J. Funct. Anal.},
   volume={253},
   date={2007},
   number={1},
   pages={220--240}
}

\bib{MR2105955}{article}{
   author={Treil, Sergei},
   title={Lower bounds in the matrix Corona theorem and the codimension one
   conjecture},
   journal={Geom. Funct. Anal.},
   volume={14},
   date={2004},
   number={5},
   pages={1118--1133}
}

\bib{MR2106344}{article}{
   author={Treil, Sergei},
   title={An operator Corona theorem},
   journal={Indiana Univ. Math. J.},
   volume={53},
   date={2004},
   number={6},
   pages={1763--1780}
}

\bib{MR1183608}{article}{
   author={Treil, Sergei},
   title={The stable rank of the algebra $H^\infty$ equals $1$},
   journal={J. Funct. Anal.},
   volume={109},
   date={1992},
   number={1},
   pages={130--154}
}

\bib{MR0981054}{article}{
   author={Treil, Sergei},
   title={Angles between co-invariant subspaces, and the operator corona
   problem. The Sz\H okefalvi-Nagy problem},
   language={Russian},
   journal={Dokl. Akad. Nauk SSSR},
   volume={302},
   date={1988},
   number={5},
   pages={1063--1068},
   translation={
      journal={Soviet Math. Dokl.},
      volume={38},
      date={1989},
      number={2},
      pages={394--399}
   }
}

\bib{MR2449054}{article}{
   author={Treil, Sergei},
   author={Wick, Brett D.},
   title={Analytic projections, corona problem and geometry of holomorphic
   vector bundles},
   journal={J. Amer. Math. Soc.},
   volume={22},
   date={2009},
   number={1},
   pages={55--76}
}

\bib{MR2158178}{article}{
   author={Treil, Sergei},
   author={Wick, Brett D.},
   title={The matrix-valued $H^p$ corona problem in the disk and polydisk},
   journal={J. Funct. Anal.},
   volume={226},
   date={2005},
   number={1},
   pages={138--172}
}

\bib{MR2057771}{article}{
   author={Trent, Tavan T.},
   title={A corona theorem for multipliers on Dirichlet space},
   journal={Integral Equations Operator Theory},
   volume={49},
   date={2004},
   number={1},
   pages={123--139}
}

\bib{MR1887635}{article}{
   author={Trent, Tavan T.},
   title={A new estimate for the vector valued corona problem},
   journal={J. Funct. Anal.},
   volume={189},
   date={2002},
   number={1},
   pages={267--282}
}

\bib{MR2213732}{article}{
   author={Trent, Tavan T.},
   author={Zhang, Xinjun},
   title={A matricial corona theorem},
   journal={Proc. Amer. Math. Soc.},
   volume={134},
   date={2006},
   number={9},
   pages={2549--2558}
}

\bib{MR2317961}{article}{
   author={Trent, Tavan T.},
   author={Zhang, Xinjun},
   title={A matricial corona theorem. II},
   journal={Proc. Amer. Math. Soc.},
   volume={135},
   date={2007},
   number={9},
   pages={2845--2854}
}

\bib{Uchiyama}{article}{
   author={Uchiyama, Akihito},
   title={Corona theorems for countably many functions and estimates for their solutions},
   status={preprint},
   date={1980}
}

\bib{MR1651405}{article}{
   author={Xiao, Jie},
   title={The $\overline\partial$-problem for multipliers of the Sobolev
   space},
   journal={Manuscripta Math.},
   volume={97},
   date={1998},
   number={2},
   pages={217--232}
}
\end{biblist}
\end{bibdiv}
\end{document}